%
%
%
\documentstyle{amsppt}
\input amsppt1.tex
\mag=\magstep1
\pagewidth{160truemm} \pageheight{235truemm}
\NoBlackBoxes
\TagsOnRight
\nologo

\define\F{\flushpar}

\define\qq{\qquad}
\define\q{\quad}
\define\SP{\smallpagebreak}
\define\MP{\medpagebreak}
\define\BP{\bigpagebreak}

\define\HH{{\Cal H}}

\define\la{{\lambda}}
\define\La{{\Lambda}}
\define\ep{{\epsilon}}
\define\de{{\delta}}
\define\lan{{\langle}}
\define\ran{{\rangle}}
\define\al{{\alpha}}
\define\be{{\beta}}
\define\sig{{\sigma}}

\define\parti{{\partial}}

\define\ellj{{\ell(j)}}

\define\RR{{\Cal R}}

\define\wtP{{\widetilde P}}
\define\whP{{\widehat P}}
\define\Mdm{{M_d^m}}

\define\TT{{\Cal T}}


\rightline{KIMS-2000-01-01}
\rightline{SP/9912244}

\vskip8pt

\topmatter
\title\nofrills \bf    Scattering Spaces and a Decomposition of Continuous Spectral Subspace\\
of $N$-body Quantum Systems
\endtitle
\author Hitoshi Kitada
\endauthor
\rightheadtext{Scattering Spaces}
\affil
Department of Mathematical Sciences,
University of Tokyo\\
Komaba, Meguro, Tokyo 153-8914, Japan\\
e-mail: kitada\@ms.u-tokyo.ac.jp,
http://kims.ms.u-tokyo.ac.jp/
\endaffil
\dedicatory
Dedicated to Professor Yoshimi Sait\B{o} on his 60th birthday\\
\\ 
{\rm December 7, 1999}\\
\enddedicatory
\abstract{
We introduce the notion of scattering space $S_b^r$ for $N$-body quantum mechanical systems, where $b$ is a cluster decomposition with $2\le |b|\le N$ and $r$ is a real number $0\le r\le 1$. Utilizing these spaces, we give a decomposition of continuous spectral subspace by $S_b^1$ for $N$-body quantum systems with long-range pair potentials $V_\alpha^L(x_\alpha)=O(|x_\al|^{-\ep})$. This is extended to a decomposition by $S_b^r$ with $0\le r\le 1$ for some long-range case. We also prove a characterization of ranges of wave operators by $S_b^0$.}
\endtopmatter

\document


\subhead\nofrills
Introduction
\endsubhead
\par
\vskip 12pt

After Derezi\'nski [D] proved the asymptotic completeness of $N$-body quantum systems with long-range pair potentials $V_\alpha^L(x_\alpha)=O(|x_\alpha|^{-\ep})$ for $\ep>\sqrt{3}-1$, Yafaev [Y] gave an example where the asymptotic completeness breaks down for $0<\ep<1/2$. It seems that the cause which prevents the asymptotic completeness in his case is some ``self-similarity" of approximate threshold eigenvectors by which property the channels corresponding to those eigenvectors cannot be separated well enough as in the usual case. The present paper is partly a reminiscence of my preprint [K3] in 1984, while apart from the problem of asymptotic completeness. After some preparations of section 1, we define in section 2 scattering spaces $S_b^r$, $2\le|b|\le N$, $0\le r\le 1$, for a given $N$-body Hamiltonian $H$. Roughly speaking $S_b^r$ consists of vectors $f$ such that in its evolution $e^{-itH}f$, the distance between clusters in $b$ evolves linearly in $t$, whereas the size of each cluster in $b$ is bounded by a constant times $t^r$. We then prove in section 4 a decomposition of the continuous spectral subspace as an orthogonal sum of the scattering spaces $S_b^1$, $2\le|b|\le N$. We extend the decomposition to the one by $S_b^r$ with general order $0\le r\le 1$ for some long-range pair potentials by utilizing the results of Derezi\'nski [D]. We also see that the orthogonal projection $E_b(r)$ onto the scattering space $S_b^r$ $(0\le r\le 1)$ is discontinuous at some $r=r_0\in(\ep,1/2)$ in the case of Yafaev's potentials [Y]. We prove in section 5 a characterization of the ranges of wave operators: $\RR(W_b^\pm)=S_b^0$ for general long-range potentials, where we adopt the wave operators $W_b^\pm$ with time-independent modifiers $J_b$ extending [IK] to $N$-body case. We assume some decay assumptions on threshold eigenvectors (Assumption 1.2) in addition to the usual assumptions on the decay of pair potentials (Assumption 1.1).

\BP

\subhead\nofrills
1. Preliminaries
\endsubhead
\par
\vskip 12pt

We consider the Schr\"odinger operator defined in $L^2(R^{\nu N})$ ($\nu\ge 1$, $N\ge 2$)
$$
H=H_0+V,\q H_0=-\sum_{i=1}^N \frac{\hbar^2}{2m_i}\frac{\parti^2}{\parti r_i^2}.\tag 1.1
$$ 
Here
$$
V=\sum_\al V_\al(x_\al),\tag 1.2
$$
where $x_\al=r_i-r_j$, $r_i=(r_{i1},\cdots,r_{i\nu})\in R^\nu$ is the position vector of the $i$-th particle, $\frac{\parti}{\parti r_i}=\left(\frac{\parti}{\parti r_{i1}},\cdots,\frac{\parti}{\parti r_{i\nu}}\right)$, $\frac{\parti^2}{\parti r_i^2}=\sum_{j=1}^\nu \frac{\parti^2}{\parti r_{ij}^2}=\Delta_{r_i}$, $m_i>0$ is the mass of the $i$-th particle, and $\al=\{i,j\}$ is a pair with $1\le i<j\le N$. Our assumption on the decay rate of the pair potentials $V_\al(x_\al)$ is as follows.

\definition{Assumption 1.1} $V_\al(x)$ ($x\in R^\nu$) is split into a sum of a real-valued $C^\infty$ function $V_\al^L(x)$ and a real-valued measurable function $V_\al^S(x)$ of $x\in R^\nu$ satisfying the following conditions: There are real numbers $\ep$ and $\ep_1$ with $0<\ep,\ep_1<1$ such that for all multi-indices $\beta$
$$
|\parti_{x}^\be V_\al^L(x)|\le C_\be\lan x\ran^{-|\be|-\ep}\tag 1.3
$$
with some constants $C_\be>0$ independent of $x\in R^\nu$, and
$$
\lan x\ran^{1+\ep_1}V_\al^S(x)(-\Delta_x+1)^{-1} \ \text{is a bounded operator in}\ L^2(R^\nu).\tag 1.4
$$
Here $\Delta_x$ is a Laplacian with respect to $x$, and $\lan x\ran$ is a $C^\infty$ function of $x$ such that
$
\lan x\ran
=|x|$ for $|x|\ge1$ and
$\ge \frac{1}{2}$ for $|x|<1$.
\enddefinition

We can adopt weaker conditions on the differentiability and decay rate for higher derivatives of the long-range part $V_\al^L(x)$, but for later convenience of exposition, we adopt this form in the present paper.

The free part $H_0$ of $H$ in (1.1) has various forms in accordance with our choice of coordinate systems. We use the so-called Jacobi coordinates. The center of mass of our $N$-particle system is
$$
X_C=\frac{m_1r_1+\cdots+m_Nr_N}{m_1+\cdots+m_N},
$$
and the Jacobi coordinates are defined by
$$
x_i=r_{i+1}-\frac{m_1r_1+\cdots+m_i r_i}{m_1+\cdots+m_i},\q i=1,2,\cdots,N-1.\tag 1.5
$$
Accordingly the corresponding canonically conjugate momentum operators are defined by
$$
P_C=\frac{\hbar}{i}\frac{\parti}{\parti X_C},\q
p_i=\frac{\hbar}{i}\frac{\parti}{\parti x_i}.
$$
Using these new $X_C,P_C,x_i,p_i$, we can rewrite $H_0$ in (1.1) as
$$
H_0={\tilde H_0}+H_C.\tag 1.6
$$
Here
$$
{\tilde H}_0=\sum_{i=1}^{N-1} \frac{1}{2\mu_i}p_i^2=-\sum_{i=1}^{N-1} \frac{\hbar^2}{2\mu_i}\Delta_{x_i},\q H_C=\frac{1}{\sum_{j=1}^N m_j}P_C^2,
$$
where $\mu_i>0$ is the reduced mass defined by the relation:
$$
\frac{1}{\mu_i}=\frac{1}{m_{i+1}}+\frac{1}{m_1+\cdots+m_i}.
$$
The new coordinates give a decomposition $L^2(R^{\nu N})=L^2(R^\nu)\otimes L^2(R^n)$ with $n=\nu(N-1)$ and in this decomposition, $H$ is decomposed
$$
H=H_C\otimes I+I\otimes {\tilde H},\q {\tilde H}={\tilde H}_0+V.
$$
$H_C$ is a Laplacian, so we consider ${\tilde H}$ in the Hilbert space $\HH=L^2(R^n)=L^2(R^{\nu(N-1)})$. We write this ${\tilde H}$ as $H$ in the followings:
$$
H=H_0+V=\sum_{i=1}^{N-1}\frac{1}{2\mu_i}p_i^2+\sum_{\al}V_\al(x_\al)=-\sum_{i=1}^{N-1}\frac{\hbar^2}{2\mu_i}\Delta_{x_i}+\sum_\al V_\al(x_\al).\tag 1.7
$$
This means that we consider the Hamiltonian $H$ in (1.1) restricted to the subspace of $R^{\nu N}$:
$$
(m_1+\cdots+m_N)X_C=m_1r_1+\cdots+m_N r_N=0.\tag 1.8
$$
We equip this subspace with the inner product:
$$
\lan x,y\ran=\sum_{i=1}^{N-1} \mu_i x_i\cdot y_i,\tag 1.9
$$
where $\cdot$ denotes the Euclidean scalar product. With respect to this inner product, the changes of variables between Jacobi coordinates in (1.5) are realized by orthogonal transformations on the space $R^n$ defined by (1.8), while $\mu_i$ and $x_i$ depend on the order of the construction of the Jacobi coordinates in (1.5). If we define velocity operator $v$ by
$$
v=(v_1,\cdots,v_{N-1})=(\mu_1^{-1}p_1,\cdots,\mu_{N-1}^{-1}p_{N-1}),
$$
we can write using the inner product above
$$
H_0=\frac{1}{2}\lan v,v\ran.\tag 1.10
$$

Next we introduce clustered Jacobi coordinate.
 Let $a=\{C_1,\cdots,C_k\}$ be a disjoint decomposition of the set $\{1,2,\cdots,N\}$: 
$C_j \ne \emptyset$
 $(j=1,2,\cdots,k)$, $\cup_{j=1}^k C_j=\{1,2,\cdots,N\}$ with $C_i\cap C_j=\emptyset$ when $i\ne j$. We denote the number of elements of a set $S$ by $|S|$. Then $|a|=k$ in the present case, and we call $a$ a cluster decomposition with $|a|$ clusters $C_1,\cdots,C_{|a|}$.
A clustered Jacobi coordinate $x=(x_a,x^a)$ associated with a cluster decomposition $a=\{C_1,\cdots,C_k\}$ is obtained by first choosing a Jacobi coordinate
$$
x^{(C_\ell)}=(x_1^{(C_\ell)},\cdots,x_{|C_\ell|-1}^{(C_\ell)})\in R^{\nu(|C_\ell|-1)}\q (\ell=1,2,\cdots,k)
$$
for the $|C_\ell|$ particles in the cluster $C_\ell$ and then by choosing an intercluster Jacobi coordinate
$$
x_a=(x_1,\cdots,x_{k-1})\in R^{\nu(k-1)}
$$
for the center of mass of the $k$ clusters $C_1,\cdots,C_k$. Then $x^a=(x^{(C_1)},\cdots,x^{(C_k)})\in R^{\nu(N-k)}$ and $x=(x_a,x^a)\in R^{\nu(N-1)}=R^n$. The corresponding canonically conjugate momentum operator is
$$\align
&p=(p_a,p^a),\q p_a=(p_1,\cdots,p_{k-1}),\q p^a=(p^{(C_1)},\cdots,p^{(C_k)}),\\
&p_i=\frac{\hbar}{i}\frac{\parti}{\parti x_i},\q
p^{(C_\ell)}=(p_1^{(C_\ell)},\cdots,p_{|C_\ell|-1}^{(C_\ell)}),\q p_i^{(C_\ell)}=\frac{\hbar}{i}\frac{\parti}{\parti x_i^{(C_\ell)}}.
\endalign
$$
Accordingly $\HH=L^2(R^n)$ is decomposed:
$$\align
&\HH=\HH_a\otimes \HH^a,\q \HH_a=L^2(R^{\nu(k-1)}_{x_a}),\\ 
&\HH^a=L^2(R^{\nu(N-k)}_{x^a})=\HH^{(C_1)}\otimes\cdots\otimes \HH^{(C_k)},\q
\HH^{(C_\ell)}=L^2(R^{\nu(|C_\ell|-1)}_{x^{(C_\ell)}}).\tag 1.11
\endalign
$$
In this coordinates system, $H_0$ in (1.7) is decomposed:
$$
\align
&H_0=T_a+H_0^a,\\
&T_a=-\sum_{\ell=1}^{k-1}\frac{\hbar^2}{2M_\ell}\Delta_{x_\ell},\tag 1.12\\
&H_0^a=\sum_{\ell=1}^k H_0^{(C_\ell)},\q H_0^{(C_\ell)}=-\sum_{i=1}^{|C_\ell|-1}\frac{\hbar^2}{2\mu_i^{(C_\ell)}}\Delta_{x_i^{(C_\ell)}},
\endalign
$$
where
$\Delta_{x_\ell}$ and $\Delta_{x_i^{(C_\ell)}}$ are $\nu$-dimensional Laplacians and $M_\ell$ and $\mu_i^{(C_\ell)}$ are the reduced masses. We introduce the inner product in the space $R^n=R^{\nu(N-1)}$ as in (1.9):
$$
\align
&\lan x,y\ran=\lan (x_a,x^a),(y_a,y^a)\ran=\lan x_a,y_a\ran+\lan x^a,y^a\ran\\
&=\sum_{\ell=1}^{k-1}M_\ell x_\ell\cdot y_\ell + \sum_{\ell=1}^k\sum_{i=1}^{|C_\ell|-1}\mu_i^{(C_\ell)}x_i^{(C_\ell)}\cdot y_i^{(C_\ell)},\tag 1.13
\endalign
$$
and velocity operator
$$
v=(v_a,v^a)=M^{-1}p=(m_a^{-1}p_a,(\mu^a)^{-1}p^a),\tag 1.14
$$
where 
$M=\left(\smallmatrix m_a&0\\ 0&\mu^a \endsmallmatrix\right)$
 is the $n=\nu(N-1)$ dimensional diagonal mass matrix whose diagonals are given by
$M_1,\cdots,M_{k-1},\mu_1^{(C_1)},\cdots,\mu_{|C_k|-1}^{(C_k)}$. Then $H_0$ is written as
$$
H_0=\frac{1}{2}\lan v,v\ran=T_a+H_0^a=\frac{1}{2}\lan v_a,v_a\ran+\frac{1}{2}\lan v^a,v^a\ran.\tag 1.15
$$

We need a notion of order in the set of cluster decompositions. A cluster decomposition $b$ is called a refinement of a cluster decomposition $a$, iff any $C_\ell\in b$ is a subset of some $D_k\in a$. When $b$ is a refinement of $a$ we denote this as $b\le a$. $b\not\leq a$ is its negation: some cluster $C_\ell\in b$ is not a subset of any $D_k\in a$. Thus for a pair $\al=\{i,j\}$, $\al\le a$ means that $\al=\{i,j\}\subset D_k$ for some $D_k\in a$, and $\al\not\leq a$ means that $\al=\{i,j\}\not\subset D_k$ for any $D_k\in a$. $b<a$ means that $b\le a$ but $b\ne a$.

We decompose the potential term $V$ in (1.7) as
$$
\sum_\al V_\al(x_\al)=V_a+I_a,\tag 1.16
$$
where
$$
\align
&V_a=\sum_{C_\ell\in a}V_{C_\ell},\\
&V_{C_\ell}=\sum_{\al\subset C_\ell}V_\al(x_\al),\tag 1.17\\
&I_a=\sum_{\al\not\leq a} V_\al(x_\al).
\endalign
$$
By definition, $V_{C_\ell}$ depends only on the variable $x^{(C_\ell)}$ inside the cluster $C_\ell$. Similarly, $V_a$ depends only on the variable $x^a=(x^{(C_1)},\cdots,x^{(C_k)})\in R^{3(N-|a|)}$, while $I_a$ depends on all components of the variable $x$.

Then $H$ in (1.7) is decomposed:
$$\align
&H=H_a + I_a=T_a\otimes I+ I\otimes H^a + I_a,\\
&H_a=H-I_a=T_a\otimes I+ I\otimes H^a,\tag 1.18\\
&H^a=H_0^a+V_a=\sum_{C_\ell\in a} H^{(C_\ell)},\q H^{(C_\ell)}=H_0^{(C_\ell)}+V_{C_\ell},
\endalign
$$
where $T_a$ is an operator in $\HH_a=L^2(R^{\nu(k-1)}_{x_a})$, $H^a$ and $H_0^a$ are operators in $\HH^a=L^2(R^{\nu(N-k)}_{x^a})$, and $H^{(C_\ell)}$ and $H_0^{(C_\ell)}$ are operators in $\HH^{(C_\ell)}=L^2(R^{\nu(|C_\ell|-1)}_{x^{(C_\ell)}})$.

We denote by $P_a$ the orthogonal projection onto the pure point spectral subspace (or eigenspace) $\HH_{pp}^a=\HH_{pp}(H^a)(\subset\HH^a)$ for $H^a$. We use the same notation $P_a$ for the obvious extention $I\otimes P_a$ to the total space $\HH$. For $|a|=N$, we set $P_a=I$. Let $M=1,2,\cdots$ and $P_a^M$ denote an $M$-dimensional partial projection of $P_a$ such that s-$\lim_{M\to\infty}P_a^M=P_a$. We define for $\ell=1,\cdots,N-1$ and an $\ell$-dimensional multi-index $M=(M_1,\cdots,M_\ell)$ ($M_j\ge 1$)
$$
{\widehat P}_\ell^M=\left(I-\sum_{|a_\ell|=\ell}P_{a_\ell}^{M_\ell}\right)\cdots\left(I-\sum_{|a_2|=2}P_{a_2}^{M_2}\right)(I-P^{M_1}).\tag 1.19
$$
(Note that for $|a|=1$, $a=\{C\}$ with $C=\{1,2,\cdots,N\}$. Thus $P^{M_1}$ is an $M_1$-dimensional partial projection into the eigenspace of $H$.)
We further define for a $|a|$-dimensional multi-index $M_a=(M_1,\cdots,M_{|a|-1},M_{|a|})=({\widehat M}_a,M_{|a|})$
$$
{\widetilde P}_a^{M_a}=P_a^{M_{|a|}}{\widehat P}_{|a|-1}^{{\widehat M}_a},\q
2\le |a| \le N.\tag 1.20
$$
Then it is clear that
$$
\sum_{2\le |a|\le N}{\widetilde P}_a^{M_a}={\widehat P}_1^{M_1}=I-P^{M_1},\tag 1.21
$$
provided that the component $M_j$ of $M_a$ depends only on the number $j$ but not on $a$. In the following we use such $M_a$'s only.

Related with those notions, we denote by $\HH_c=\HH_c(H)$ the orthogonal complement $\HH_{pp}(H)^\perp$ of the eigenspace $\HH_{pp}=\HH_{pp}(H)$ for the total Hamiltonian $H$. Namely $\HH_c(H)$ is the continuous spectral subspace for $H$. We note that $\HH_c(H)=(I-P_a)\HH$ for a unique $a$ with $|a|=1$, and that for $f\in\HH$, $(I-P^{M_1})f\to (I-P_a)f\in \HH_c(H)$ as $M_1\to\infty$. We use freely the notations of functional analysis for selfadjoint operators, e.g. $E_H(\Delta)$ is the spectral measure for $H$.

To state a theorem due to Enss [E2], we introduce an assumption:

\definition{Assumption 1.2} For any cluster decomposition $a$ with $2\le |a|\le N-1$ and any integer $M=1,2,\cdots$,
$$
\Vert |x^a|^2P_a^M\Vert<\infty.\tag 1.22
$$
\enddefinition

This assumption is concerned with the decay rate of eigenvectors of subsystem Hamiltonians. Since it is known that non-threshold eigenvectors decay exponentially (see Froese and Herbst [FH]), this assumption is the one about threshold eigenvectors.

Let $v_a$, as above, denote the velocity operator between the clusters in $a$. It is expressed as $v_a=m_a^{-1}p_a$ for some $\nu(|a|-1)$-dimensional diagonal mass matrix $m_a$. Then we can state the theorem.

\proclaim{Theorem 1.3([E2])}
Let $N\ge 2$ and let $H$ be the Hamiltonian $H$ in (1.7) or (1.18) for an $N$-body quantum-mechanical system. Let Assumptions 1.1 and 1.2 be satisfied.
Let $f\in \HH$. Then there exist a sequence $t_m\to\pm\infty$ (as $m\to\pm\infty$) and a sequence $M_a^m$ of multi-indices whose components all tend to $\infty$ as $m\to\pm\infty$ such that for all cluster decompositions $a$ with $2\le |a|\le N$, for all $\varphi\in C_0^\infty(R_{x_a}^{\nu(|a|-1)})$, $R>0$, and $\alpha=\{i,j\}\not\leq a$
$$\align
&\left\Vert\frac{|x^a|^2}{t_m^2}{\widetilde P}_a^{M_a^m}e^{-it_mH/\hbar}f\right\Vert \to 0\tag 1.23\\
&\Vert F(|x_\alpha|<R){\widetilde P}_a^{M_a^m}e^{-it_mH/\hbar}f\Vert \to 0\tag 1.24\\
&\Vert (\varphi(x_a/t_m)-\varphi(v_a)){\widetilde P}_a^{M_a^m}e^{-it_mH/\hbar}f\Vert \to 0\tag 1.25
\endalign
$$
as $m\to\pm\infty$. Here $F(S)$ is the characteristic function of the set defined by the condition $S$.
\endproclaim

We denote the sum of the sets of thresholds and eigenvalues of $H$ by $\TT$:
$$
\TT=\bigcup_{1\le |a|\le N}\sig_p(H^a)={\tilde \TT}\cup\sig_p(H),\q {\tilde \TT}=\bigcup_{2\le |a|\le N}\sig_p(H^a)\tag 1.26
$$
where 
$$
\sig_p(H^a)=\{\tau_1+\cdots+\tau_{|a|}\ |\ \tau_\ell\in\sig_p(H^{(C_\ell)})\ (C_\ell\in a)\}\tag 1.27
$$
is the set of eigenvalues of a subsystem Hamiltonian $H^a=\sum_{C_\ell\in a}H^{(C_\ell)}$. For $|a|=N$ we define $\sig_p(H^a)=\{0\}$. Similarly $\TT_a$ and ${\tilde \TT}_a$ are defined:
$$
\TT_a=\bigcup_{b\le a}\sig_p(H^b)={\tilde \TT}_a\cup\sig_p(H^a),\q {\tilde \TT}_a=\bigcup_{b<a}\sig_p(H^b).\tag 1.28
$$
It is known (Froese and Herbst [FH]) that these sets are subsets of $(-\infty,0]$. Further these sets form bounded, closed and countable subsets of $R^1$, and $\sig_p(H^a)$ accumulates only at ${\tilde \TT}_a$ (see Cycon {\it et al}. [C]).

We use the notation $\Delta\Subset \Delta'$ for Borel sets $\Delta, \Delta'\subset R^k$ to mean that the closure $\bar{\Delta}$ of $\Delta$ is compact in $R^k$ and is a subset of the interior of $\Delta'$.
\BP

\vskip 12pt

\BP

\subhead\nofrills
2. Scattering Spaces
\endsubhead
\par
\vskip 12pt

In the following we consider the case $t\to\infty$ only. The other case $t\to-\infty$ is treated similarly. We also choose a unit system such that $\hbar=1$. We use the notation $f(t)\sim g(t)$ as $t\to\infty$ to mean that $\Vert f(t)-g(t)\Vert\to 0$ as $t\to\infty$ for $\HH$-valued functions $f(t)$ and $g(t)$ of $t>1$.

\definition{Definition 2.1} Let real numbers $r, \sig, \de$ and a cluster decomposition $b$ satisfy $0\le r\le 1$, $\sig, \de >0$ and $2\le |b|\le N$. 
\BP

\F
i) Let $\Delta\Subset R^1-\TT$ be a closed set.
We define $S_b^{r\sig\de}(\Delta)$ for $0<r\le 1$ by
$$
S_b^{r\sig\de}(\Delta)=\{f\in E_H(\Delta)\HH \ |\ 
 e^{-itH}f\sim \prod_{\al\not\leq b}F(|x_\al|\ge \sig t)F(|x^b|\le\de t^r)e^{-itH}f\ \text{as}\ t\to\infty\}.\tag 2.1
$$
For $r=0$ we define $S_b^{0\sig}(\Delta)$ by
$$\align
S_b^{0\sig}(\Delta)=\{&f\in E_H(\Delta)\HH \ |\ \\
&\lim_{R\to\infty}\limsup_{t\to\infty}\Bigl\Vert e^{-itH}f - \prod_{\al\not\leq b}F(|x_\al|\ge \sig t)F(|x^b|\le R)e^{-itH}f\Bigr\Vert=0 \}.\tag 2.2
\endalign
$$
We then define the localized scattering space $S_b^r(\Delta)$ of order $r\in(0,1]$ for $H$ as the closure of
$$
\align
\bigcup_{\sig>0}\bigcap_{\de>0}S_b^{r\sig\de}(\Delta)
=\{&f\in  E_H(\Delta)\HH \ |\ \exists \sig>0, \forall \de>0:\\
&\ e^{-itH}f\sim \prod_{\al\not\leq b}F(|x_\al|\ge \sig t)F(|x^b|\le\de t^r)e^{-itH}f\ \text{as}\ t\to\infty\}.\tag 2.3
\endalign
$$
$S_b^0(\Delta)$ is defined as the closure of
$$\align
\bigcup_{\sig>0}S_b^{0\sig}(\Delta)=\{&f\in E_H(\Delta)\HH\ |\ \exists \sig>0:\\
&\lim_{R\to\infty}\limsup_{t\to\infty}\Bigl\Vert e^{-itH}f - \prod_{\al\not\leq b}F(|x_\al|\ge \sig t)F(|x^b|\le R)e^{-itH}f\Bigr\Vert=0 \}.\tag 2.4
\endalign
$$

\F
ii) We define the scattering space $S_b^r$ of order $r\in [0,1]$ for $H$ as the closure of
$$
\bigcup_{\Delta\Subset R^1-\TT}S_b^r(\Delta).\tag 2.5
$$
\enddefinition
We note that $S_b^{r\sig\de}(\Delta)$, $S_b^{0\sig}(\Delta)$, $S_b^{r}(\Delta)$ and $S_b^{r}$ define closed subspaces of $E_H(\Delta)\HH$ and $\HH_c(H)$, respectively. 
\BP

\proclaim{Proposition 2.2} Let $\Delta\Subset R^1-\TT$ and $f\in S_b^{r\sig\de}(\Delta)$ for $0< r\le 1$ or $f\in S_b^{0\sig}(\Delta)$ for $r=0$ with $\sig,\de>0$ and $2\le |b|\le N$. Then the following limit relations hold:
\MP

\F
i) Let $\al\not\leq b$. Then for $0<r\le 1$ we have when $t\to\infty$
$$
F(|x_\al|<\sig t)F(|x^b|\le\de t^r)e^{-itH}f\to 0.\tag 2.6
$$
For $r=0$ we have
$$
\lim_{R\to\infty}\limsup_{t\to\infty}\left\Vert F(|x_\al|< \sig t)F(|x^b|\le R)e^{-itH}f\right\Vert=0.\tag 2.7
$$
\F
ii) For $0<r\le 1$ we have when $t\to\infty$
$$
F(|x^b|> \de t^r)e^{-itH}f\to 0.\tag 2.8
$$
For $r=0$
$$
\lim_{R\to\infty}\limsup_{t\to\infty}\left\Vert F(|x^b|> R)e^{-itH}f\right\Vert=0.\tag 2.9
$$
\F
iii) There exists a sequence $t_m\to\infty$ as $m\to\infty$ depending on $f\in S_b^{r\sig\de}(\Delta)$ or $f\in S_b^{0\sig}(\Delta)$ such that
$$
\left\Vert \left(\varphi\left({x_b}/{t_m}\right)-\varphi(v_b)\right)e^{-it_mH}f\right\Vert\to 0\q \text{as}\q m\to\infty \tag 2.10
$$
for any function $\varphi\in C_0^\infty(R^{\nu(|b|-1)}_{x_b})$.
\endproclaim
\demo{Proof}
i) and ii) are clear from the definition of $S_b^{r\sig\de}(\Delta)$ or $S_b^{0\sig}(\Delta)$. We prove iii). 
Since $f\in E_H(\Delta)\HH\subset H_c(H)$, we have by (1.21), Theorem 1.3 and $f\in S_b^{r\sig\de}(\Delta)$ (or $f\in S_b^{0\sig}(\Delta)$)
$$
e^{-it_mH}f\sim \sum_{d\le b}{\widetilde P}_d^{M_d^m}e^{-it_mH}f\tag 2.11
$$
along some sequence $t_m\to\infty$ depending on $f$. On each state on the right-hand side (RHS) of (2.11), (1.25) with $a$ replaced by $d$ holds. By the restriction $d\le b$ in the sum of the RHS of (2.11), we obtain (2.10).\ $\square$

\enddemo
\BP

The following propositions are obvious by definition.

\proclaim{Proposition 2.3} Let $2\le |b|\le N$. If $1\ge r'\ge r> 0$, $\sig\ge \sig'>0$ and $\de'\ge \de>0$ and $\Delta\Subset R^1-\TT$, then $S_b^{0\sig}(\Delta)\subset S_b^{0\sig'}(\Delta)$, $S_b^{0\sig}(\Delta)\subset S_b^{r\sig\de}(\Delta)\subset S_b^{r'\sig'\de'}(\Delta)$, $S_b^0(\Delta)\subset S_b^r(\Delta)\subset S_b^{r'}(\Delta)$, $S_b^0(\Delta)\subset S_b^r(\Delta)\subset S_b^{r}$, and $S_b^0\subset S_b^r\subset S_b^{r'}$.
\endproclaim

\proclaim{Proposition 2.4} Let $b$ and $b'$ be different cluster decompositions: $b\ne b'$. Then for any $0\le r \le 1$, $S_b^r$ and $S_{b'}^r$ are orthogonal mutually: $S_b^r\perp S_{b'}^r$.
\endproclaim

\vskip 12pt

\BP

\subhead\nofrills
3. A Partition of Unity
\endsubhead
\par
\vskip 12pt

To state a proposition that will play a fundamental role in our decomposition of continuous spectral subspace by $S_b^1$, we prepare some notations. Let $b$ be a cluster decomposition with $2\le|b|\le N$. For any two clusters $C_1$ and $C_2$ in $b$, we define a vector $z_{b1}$ that connects the two centers of mass of the clusters $C_1$ and $C_2$. The number of such vectors when we move over all pairs of clusters in $b$ is $k_b=\left(\smallmatrix |b|\\ 2\endsmallmatrix\right)$ in total. We denote these vectors by $z_{b1},z_{b2},\cdots,z_{bk_b}$.

 Let $z_{bk}$ $(1\le k\le k_b)$ connect two clusters $C_\ell$ and $C_m$ in $b$ $(\ell\ne m)$. Then for any pair $\al=\{i,j\}$ with $i\in C_\ell$ and $j\in C_m$, the vector $x_\al=x_{ij}$ is expressed like $(z_{bk},x^{(C_\ell)},x^{(C_m)})\in R^{\nu+\nu(|C_\ell|-1)+\nu(|C_m|-1)}$, where $x^{(C_\ell)}(\in R^{\nu(|C_\ell|-1)})$ and $x^{(C_m)}(\in R^{\nu(|C_m|-1)})$ are the positions of the particles $i$ and $j$ in $C_\ell$ and $C_m$, respectively. The vector expression $x_\al=(z_{bk},x^{(C_\ell)},x^{(C_m)})$ is in the space $R^{\nu+\nu(|C_\ell|-1)+\nu(|C_m|-1)}$. If we express it in the larger space $R^{\nu}_{z_{bk}}\times R^{\nu(N-|b|)}_{x^b}$, it would be $x_\al=(z_{bk},x^b)$, and $|x_\al|^2=|z_{bk}|^2+|x^b|^2$.  Thus if $|z_{bk}|^2$ is sufficiently large compared to $|x^b|^2\ge|x^{(C_\ell)}|^2+|x^{(C_m)}|^2$, e.g. if $|z_{bk}|^2> \rho>0$ and $|x^b|^2< \theta$ with $\rho\gg \theta>0$ (which means that $\rho/\theta$ is sufficiently large), then $|x_\al|^2> \rho/2$ for all $\al\not\leq b$. 

Next if $c<b$ and $|c|=|b|+1$, then just one cluster, say $C_\ell\in b$, is decomposed into two clusters $C'_\ell$ and $C''_\ell$ in $c$, and other clusters in $b$ remain the same in the finer cluster decomposition $c$. In this case, we can choose just one vector $z_{ck}$ $(1\le k\le k_c)$ that connects clusters $C'_\ell$ and $C''_\ell$ in $c$, and we can express $x^b=(z_{ck},x^c)$. The norm of this vector is written as 
$$
|x^b|^2=|z_{ck}|^2+|x^c|^2.\tag 3.1
$$
Similarly the norm of $x=(x_b,x^b)$ is written as
$$
|x|^2=|x_b|^2+|x^b|^2.\tag 3.2
$$
 We recall that norm is defined, as usual, from the inner product defined by (1.13) which changes in accordance with the cluster decomposition used in each context. E.g., in (3.1), the left-hand side (LHS) is defined by using (1.13) for the cluster decomposition $b$ and the RHS is by using (1.13) for $c$.

With these preparations, we state the following lemma, which is partially a repetition of [K1, Lemma 2.1]. We define subsets $T_b(\rho,\theta)$ and ${\tilde T}_b(\rho,\theta)$ of $R^n=R^{\nu(N-1)}$ for cluster decompositions $b$ with $2\le |b|\le N$ and real numbers $\rho,\theta$ with $1>\rho,\theta>0$:
$$\align
&T_b(\rho,\theta)=\left(\bigcap_{k=1}^{k_b}\{x\ |\ |z_{bk}|^2>\rho|x|^2\}\right)\cap\{ x\ |\ |x_b|^2>(1-\theta)|x|^2\},\tag 3.3\\
&{\tilde T}_b(\rho,\theta)=\left(\bigcap_{k=1}^{k_b}\{x\ |\ |z_{bk}|^2>\rho\}\right)\cap\{ x\ |\ |x_b|^2>1-\theta\}.\tag 3.4
\endalign
$$
Subsets $S$ and $S_{\theta}$ $(\theta>0)$ of $R^n=R^{\nu(N-1)}$ are defined by
$$\align
S&=\{x\ |\ |x|^2\ge 1\},\\
S_\theta&=\{x\ |\ 1+\theta\ge |x|^2\ge 1\}.
\endalign
$$

\proclaim{Lemma 3.1} 
Suppose that constants $1\ge\theta_1>\rho_j>\theta_j> \rho_N>0$ satisfy $\theta_{j-1}\ge \theta_j+\rho_j$ for $j=2,3,\cdots,N-1$. 
Then the followings hold:
\BP

\F
i)
$$
S\subset \bigcup_{2\le |b|\le N}T_b(\rho_{|b|},\theta_{|b|}).\tag 3.5
$$
\MP

\F
ii)
Let $\gamma_j>1$ $(j=1,2)$ satisfy
$$
\gamma_1\gamma_2< r_0:=\min_{2\le j\le N-1}\{{\rho_j/\theta_j}\}.\tag 3.6
$$
If $b\not\le c$ with $|b|\ge |c|$, then
$$
T_b(\gamma_1^{-1}\rho_{|b|},\gamma_2\theta_{|b|})\cap
T_c(\gamma_1^{-1}\rho_{|c|},\gamma_2\theta_{|c|})=\emptyset.\tag 3.7
$$
\MP

\F
iii)
For $\gamma>1$ and $2\le|b|\le N$
$$\align
T_b(\rho_{|b|},\theta_{|b|})\cap S_{\theta_{N-1}}&\subset
{\tilde T}_b(\rho_{|b|},\theta_{|b|})\cap S_{\theta_{N-1}}\\
&\Subset {\tilde T}_b(\gamma^{-1}\rho_{|b|},\gamma\theta_{|b|})\cap S_{\theta_{N-1}}\\
&\subset T_b({\gamma'_1}^{-1}\rho_{|b|},\gamma'_2\theta_{|b|})\cap S_{\theta_{N-1}},
\tag 3.8
\endalign
$$
where
$$
\gamma'_1=\gamma(1+\theta_{N-1}),\q \gamma'_2=(1+\gamma)(1+\theta_{N-1})^{-1}.\tag 3.9
$$
\MP

\F
iv)
If $\frac{2\gamma'_1\gamma'_2}{2-\gamma'_1} <r_0$, then for $2\le|b|\le N$
$$
{T}_b({\gamma'_1}^{-1}\rho_{|b|},\gamma'_2\theta_{|b|})\subset \{ x\ |\ |x_\al|^2> \rho_{|b|}|x|^2/2 \ \text{for all}\ \al\not\leq b\}.\tag 3.10
$$
\MP

\F
v)
If $\gamma(1+\gamma)<r_0$ and $b\not\le c$ with $|b|\ge |c|$, then
$$
T_b({\gamma'_1}^{-1}\rho_{|b|},\gamma'_2\theta_{|b|})\cap
T_c({\gamma'_1}^{-1}\rho_{|c|},\gamma'_2\theta_{|c|})=\emptyset.\tag 3.11
$$
\endproclaim
\demo{Proof} To prove (3.5), suppose that $|x|^2\ge 1$ and $x$ does not belong to the set
$$
A=\bigcup_{2\le |b|\le N-1}\left[\left(\bigcap_{k=1}^{k_b}\{x\ |\ |z_{bk}|^2>\rho_{|b|}|x|^2\}\right)\cap\{ x\ |\ |x_b|^2>(1-\theta_{|b|})|x|^2\}\right].
$$
Under this assumption, we prove $|x_\al|^2>\rho_N|x|^2$ for all pairs $\al=\{i,j\}$. (Note that $z_{bk}$ for $|b|=N$ equals some $x_\al$.) Let $|b|=2$ and write $x=(z_{b1},x^b)$. Then by (3.1), $1\le|x|^2=|z_{b1}|^2+|x^b|^2$. Since $x$ belongs to the complement $A^c$ of the set $A$, we have $|z_{b1}|^2\le \rho_{|b|}|x|^2$ or $|x_b|^2\le(1-\theta_{|b|})|x|^2$. If $|z_{b1}|^2\le \rho_{|b|}|x|^2$, then $|x^b|^2=|x|^2-|z_{b1}|^2\ge(1-\rho_{|b|})|x|^2\ge(\theta_1-\rho_{|b|})|x|^2\ge \theta_{|b|}|x|^2$ by $\theta_{j-1}\ge \theta_j+\rho_j$. Thus $|x_b|^2=|x|^2-|x^b|^2\le(1-\theta_{|b|})|x|^2$ for all $b$ with $|b|=2$.

Next let $|c|=3$ and assume $|x_c|^2>(1-\theta_{|c|})|x|^2$. Then by $x\in A^c$, we can choose $z_{ck}$ with $1\le k\le k_c$ such that $|z_{ck}|^2\le \rho_{|c|}|x|^2$. Let $C_\ell$ and $C_m$ be two clusters in $c$ connected by $z_{ck}$, and let $b$ be the cluster decomposition obtained by combining $C_\ell$ and $C_m$ into one cluster with retaining other clusters of $c$ in $b$. Then $|b|=2$, $x^b=(z_{ck},x^c)$, and $|x^b|^2=|z_{ck}|^2+|x^c|^2$. Thus $|x_b|^2=|x|^2-|x^b|^2=|x|^2-|z_{ck}|^2-|x^c|^2=|x_c|^2-|z_{ck}|^2>(1-\theta_{|c|}-\rho_{|c|})|x|^2\ge(1-\theta_{|b|})|x|^2$, which contradicts the result of the previous step. Thus $|x_c|^2\le(1- \theta_{|c|})|x|^2$ for all $c$ with $|c|=3$. 

Repeating this procedure, we finally arrive at $|x_d|^2\le (1-\theta_{|d|})|x|^2$, thus $|x^d|^2=|x|^2-|x_d|^2\ge \theta_{|d|}|x|^2> \rho_N|x|^2$ for all $d$ with $|d|=N-1$. Namely $|x_\al|^2>\rho_N|x|^2$ for all pairs $\al=\{i,j\}$.
The proof of (3.5) is complete.

We next prove (3.7). By $b\not\le c$, we can take a pair $\al=\{i,j\}$ and clusters $C_\ell, C_m\in c$ such that $\al\le b$, $i\in C_\ell$, $j\in C_m$, and $\ell\ne m$. Then we can write $x_\al=(z_{ck},x^c)$ for some $1\le k\le k_c$. Thus if there is $x\in T_b(\gamma_1^{-1}\rho_{|b|},\gamma_2\theta_{|b|})\cap T_c(\gamma_1^{-1}\rho_{|c|},\gamma_2\theta_{|c|})$, then 
$$
\gamma_2\theta_{|b|}|x|^2>|x^b|^2\ge|x_\al|^2=|z_{bk}|^2+|x^c|^2\ge |z_{ck}|^2>\gamma_1^{-1}\rho_{|c|}|x|^2.\tag 3.12
$$
 But since $|b|\ge |c|$, we have $\rho_{|c|}> \gamma_1\gamma_2\theta_{|b|}$ when $|b|=|c|$ by (3.6), and $\rho_{|c|}>\gamma_1\gamma_2\theta_{|c|}\ge \gamma_1\gamma_2(\theta_{|b|}+\rho_{|b|})>\gamma_1\gamma_2\theta_{|b|}$ when $|c|<|b|$ by $\theta_{j-1}\ge \theta_j+\rho_j$, which both contradict the inequality (3.12). This completes the proof of (3.7).

(3.8) follows by a simple calculation from the inequality $|x|^2(1+\theta_{N-1})^{-1}\le 1$ that holds on $S_{\theta_{N-1}}$. (3.10) follows from the relation $|x_\al|^2=|z_{bk}|^2+|x^b|^2$ stated before the lemma, and (3.11) from $\gamma'_1\gamma'_2=\gamma(1+\gamma)$ and ii).\ $\square$

\enddemo
\BP

In the followings we fix constants  $\gamma>1$ and $1\ge \theta_1>\rho_j>\theta_j>\rho_N>0$ such that
$$\align
&\theta_{j-1}\ge \theta_j+\rho_j\q(j=2,3,\cdots,N-1),\tag 3.13\\
&\max\left\{\gamma(1+\gamma),\frac{2\gamma'_1\gamma'_2}{2-\gamma'_1}\right\}<r_0=\min_{2\le j\le N-1}\{\rho_j/\theta_j\},\tag 3.14
\endalign
$$
where $\gamma'_j$ $(j=1,2)$ are defined by (3.9).

Let $\rho(\la)\in C^\infty(R^1)$ be such that $0\le \rho(\la)\le 1$, $\rho(\la)=1$ $(\la\le -1)$, $\rho(\la)=0$ $(\la\ge 0)$, and $\rho'(\la)\le 0$. Then we define  functions $\phi_\sigma(\la<\tau)$ and $\phi_\sigma(\la>\tau)$ of $\la\in R^1$ by
$$\align
&\phi_\sig(\la<\tau)=\rho((\la-(\tau+\sig))/\sig),\\
&\phi_\sig(\la>\tau)=1-\phi_\sig(\la<\tau-\sig)
\tag 3.15
\endalign
$$
for constants $\sig>0,\tau\in R^1$.
We note that $\phi_\sig(\la<\tau)$ and $\phi_\sigma(\la>\tau)$ satisfy
$$\align
&\phi_\sig(\la<\tau)=
\cases
1 \q&(\la\le\tau)\\
0\q& (\la\ge\tau+\sig)
\endcases
\tag 3.16\\
&\phi_\sigma(\la>\tau)=
\cases
0 \q &(\la\le\tau-\sig)\\
1 \q& (\la\ge\tau)
\endcases
\tag 3.17\\
&\phi'_\sig(\la<\tau)=\frac{d}{d\la}\phi_\sig(\la<\tau)\le 0,\\
&\phi'_\sig(\la>\tau)\ge 0.
\endalign
$$

We define for a cluster decomposition $b$ with $2\le |b|\le N$ 
$$
\varphi_b(x_b)=\prod_{k=1}^{k_b}\phi_\sig(|z_{bk}|^2>\rho_{|b|})\phi_\sig(|x_b|^2>1-\theta_{|b|}),\tag 3.18
$$
where $\sig>0$ is fixed as
$$
0<\sig<\min_{2\le j\le N-1}\{(1-\gamma^{-1})\rho_N,(1-\gamma^{-1})\rho_j,(\gamma-1)\theta_j\}.\tag 3.19
$$
Then $\varphi_b(x_b)$ satisfies for $x\in S_{\theta_{N-1}}$
$$ 
\varphi_b(x_b)=
\cases
1\q \text{for}\q x\in {\tilde T}_b(\rho_{|b|},\theta_{|b|}),\\
0\q \text{for}\q x\not\in {\tilde T}_b(\gamma^{-1}\rho_{|b|},\gamma\theta_{|b|}).
\endcases\tag 3.20
$$
We set for $|b|=k$ $(k=2,3,\cdots,N)$
$$
J_b(x)=\varphi_b(x_b)
\left(1-\sum_{|b_{k-1}|=k-1}\varphi_{b_{k-1}}(x_{b_{k-1}})\right)
\cdots
\left(1-\sum_{|b_{2}|=2}\varphi_{b_{2}}(x_{b_{2}})\right).\tag 3.21
$$
By v) and iii) of Lemma 3.1 and (3.20), the sums on the RHS remain only in the case $b<b_j$ for $j=k-1,\cdots,2$ and $x\in S_{\theta_{N-1}}$:
$$
J_b(x)=\varphi_b(x_b)
\left(1-\sum_{\Sb |b_{k-1}|=k-1\\ b<b_{k-1}\endSb}\varphi_{b_{k-1}}(x_{b_{k-1}})\right)
\cdots
\left(1-\sum_{\Sb |b_{2}|=2\\ b<b_2\endSb}\varphi_{b_{2}}(x_{b_{2}})\right).\tag 3.22
$$
Thus $J_b(x)$ is a function of the variable $x_b$ only:
$$
J_b(x)=J_b(x_b)\q\text{when}\q x=(x_b,x^b)\in S_{\theta_{N-1}}.\tag 3.23
$$
We also note that the supports of $\varphi_{b_j}$ in each sum on the RHS of 
(3.21) are disjoint mutually in $S_{\theta_{N-1}}$ by iii) and v) of Lemma 3.1.
By (3.5) and (3.8) of lemma 3.1, and the definition (3.18)-(3.21) of $J_b(x_b)$, we therefore have
$$
\sum_{2\le|b|\le N}J_b(x_b)=1\q\text{on}\q S_{\theta_{N-1}}.
$$

We have constructed a partition of unity on $S_{\theta_{N-1}}$:

\proclaim{Proposition 3.2} Let real numbers $1\ge\theta_1>\rho_j>\theta_j>\rho_N>0$ satisfy $\theta_{j-1}\ge \theta_j+\rho_j$ for $j=2,3,\cdots,N-1$. 
Assume that (3.14) hold and let $J_b(x_b)$ be defined by (3.18)-(3.22). Then we have 
$$
\sum_{2\le|b|\le N}J_b(x_b)=1\q\text{on}\q S_{\theta_{N-1}}.\tag 3.24
$$
$J_b(x_b)$ is a $C^\infty$ function of $x_b$ and satisfies $0\le J_b(x_b)\le 1$. Further on {\rm supp}\hskip2pt$J_b \cap S_{\theta_{N-1}}$ we have
$$
|x_\al|^2>\rho_{|b|}|x|^2/2\tag 3.25
$$
for any pair $\al\not\le b$, and
$$
\sup_{x\in R^n,2\le|b|\le N}|\nabla_{x_b} J_b(x_b)|<\infty\tag 3.26
$$
for each fixed $\sig>0$ in (3.18)-(3.19).
\endproclaim
\demo{Proof} We have only to see (3.25) and (3.26). But (3.25) is clear by (3.8), (3.10), (3.14), (3.20) and (3.21), and (3.26) follows from (3.15), (3.18) and (3.22).\ $\square$
\enddemo

\vskip 12pt

\BP

\subhead\nofrills
4. A Decomposition of Continuous Spectral Subspace
\endsubhead
\par
\vskip 12pt

The following theorem gives a decomposition of $\HH_c(H)$ by scattering spaces $S_b^1$ $(2\le|b|\le N)$. 

\proclaim{Theorem 4.1} Let Assumptions 1.1 and 1.2 be satisfied. Then we have
$$
\HH_c(H)=\bigoplus_{2\le|b|\le N}S_b^1.\tag 4.1
$$
\endproclaim
\demo{Proof} Since the set
$$
\bigcup_{\Delta\Subset R^1-\TT}E_H(\Delta)\HH
$$
is dense in $\HH_c(H)$, and $S_b^1$ $(2\le|b|\le N)$ are closed and mutually orthogonal, it suffices to prove that any $\Phi(H)f$ with $\Phi\in C_0^\infty(R^1-\TT)$ and $f\in\HH$ can be decomposed as a sum of the elements $f_b^1$ in $S_b^1$: $\Phi(H)f=\sum_{2\le|b|\le N}f_b^1$.

We divide the proof into two steps. In the first step I), we prove existence of certain time limits. In the second step II), we prove existence of some ``boundary values" of those limits, and conclude the proof of decomposition (4.1).
\BP

\F
I) Existence of some time limits:
\MP

We decompose $\Phi(H)f$ as a finite sum: $\Phi(H)f=\sum_{j_0}^{\text{finite}}\psi_{j_0}(H)f$, where $\psi_{j_0}\in C_0^\infty(R^1-\TT)$. In the step I), we will prove the existence of the limit
$$
\lim_{t\to\infty}\sum_{\ell=1}^L e^{itH}
G_{b,\la_\ell}(t)^*J_b(v_b/r_\ell)G_{b,\la_\ell}(t)
e^{-itH}\psi_{j_0}(H)f,\tag 4.2
$$
under the assumption that supp\hskip2pt$\psi_{j_0}\subset \tilde\Delta\Subset\Delta$ for some intervals $\tilde\Delta\Subset\Delta\Subset R^1-\TT$ with $E\in \Delta$ and diam\hskip2pt$\Delta<d(E)$, where $d(E)>0$ is some small constant depending on $E\in R^1-\TT$ and diam\hskip2pt$S$ denotes the diameter of a set $S\subset R^1$.
The relevant factors in (4.2) will be defined in the course of the proof.
We will write $f$ for $\psi_{j_0}(H)f$ in the followings.

We take $\psi\in C_0^\infty(R^1)$ such that $\psi(\la)=1$ for $\la\in \tilde\Delta$ and supp\hskip2pt$\psi\subset{\Delta}$ for the intervals ${\tilde \Delta}\Subset\Delta$ above. Then $f=\psi(H)f=E_H(\Delta)f \in E_H(\Delta)\HH\subset \HH_c(H)$ and $e^{-itH}f=\psi(H)e^{-itH}f$. Thus
 we can use the decomposition (1.21) for the sequences $t_m$ and $M_b^m$ in Theorem 1.3:
$$
e^{-it_mH}f=\psi(H)e^{-it_mH}f=\psi(H)\sum_{2\le|d|\le N}{\widetilde P}_d^{M_d^m}e^{-it_mH}f.\tag 4.3
$$
By Theorem 1.3-(1.24)
$$
\psi(H){\widetilde P}_d^{M_d^m}e^{-it_mH}f\sim \psi(H_d){\widetilde P}_d^{M_d^m}e^{-it_mH}f\tag 4.4
$$
as $m\to\infty$. Since
$$
{\widetilde P}_d^{M_d^m}=P_d^{M_{|d|}^m}{\widehat P}_{|d|-1}^{{\widehat M}_d^m},\q
P_d^{M_{|d|}^m}=\sum_{j=1}^{M_{|d|}^m}P_{d,E_j},\tag 4.5
$$
where $P_{d,E_j}$ is one dimensional eigenprojection for $H^d$ with eigenvalue $E_j$, the RHS of (4.4) equals
$$
\sum_{j=1}^{M_{|d|}^m}\psi(T_d+E_j)P_{d,E_j}{\widehat P}_{|d|-1}^{{\widehat M}_d^m}e^{-it_mH}f.\tag 4.6
$$
By supp\hskip2pt$\psi\subset{\Delta}\Subset R^1-\TT$, $E_j\in \TT$, and $T_d\ge0$, we can take constants $\Lambda_d>\lambda_d>0$ independent of $j=1,2,\cdots$ such that $\Lambda_d\ge T_d\ge \la_d$ if $\psi(T_d+E_j)\ne 0$. Set $\Lambda_0=\max_{d}\Lambda_d>\la_0=\min_{d}\la_d>0$ and
$$
\Sigma(E)=\{E-\la\ |\ \la\in\TT, E\ge \la\}.\tag 4.7
$$
 Note that we can take $\Lambda_0>\la_0>0$ so that 
$$
\Sigma(E)\Subset(\la_0,\La_0) \subset (0,\infty).\tag 4.8
$$
Let $\Psi\in C_0^\infty(R^1)$ satisfy $\Psi(\la)=1$ for $\la\in [\la_0,\Lambda_0]$ and supp\hskip2pt$\Psi\subset[\la_0-\kappa,\Lambda_0+\kappa]$ for some small constant $\kappa>0$ such that the set $[\la'_0,\la_0]\cup[\Lambda_0,\Lambda'_0]$ is bounded away from $\Sigma(E)$, where $\la'_0=\la_0-2\kappa>0$ and $\Lambda'_0=\Lambda_0+2\kappa$. Then the RHS of (4.6) equals for any $m=1,2,\cdots$
$$
\Psi^2(T_d)\psi(H_d){\widetilde P}_d^{M_d^m}e^{-it_mH}f
.\tag 4.9
$$
On the other hand, by Theorem 1.3-(1.23) and (1.25), we have
$$
\frac{|x^d|^2}{t_m^2}\psi(H_d){\widetilde P}_d^{M_d^m}e^{-it_mH}f\sim 0\tag 4.10
$$
and
$$
\Psi^2(T_d)\psi(H_d){\widetilde P}_d^{M_d^m}e^{-it_mH}f\sim
\Psi^2(|x_d|^2/(2t_m^2))\psi(H_d){\widetilde P}_d^{M_d^m}e^{-it_mH}f\tag 4.11
$$
as $m\to\infty$, where to see (4.10) we used (1.23) and $i[H^d,|x^d|^2/t^2]=i[H_0^d,|x^d|^2/t^2]=2A^d/t^2$ where $A^d=(x^d\cdot p^d+p^d\cdot x^d)/2$, and to see (4.11) the fact that $|x_d|^2/t_m^2$ and $H_d$ commute asymptotically as $m\to\infty$ by (1.25). Thus by $|x|^2=|x_d|^2+|x^d|^2$ we have
$$
\Psi^2(T_d)\psi(H_d){\widetilde P}_d^{M_d^m}e^{-it_mH}f\sim
\Psi^2(|x|^2/(2t_m^2))\psi(H_d){\widetilde P}_d^{M_d^m}e^{-it_mH}f\tag 4.12
$$
as $m\to\infty$. From (4.3)-(4.4), (4.6), (4.9) and (4.12), we obtain
$$
e^{-it_mH}f\sim
\Psi^2(|x|^2/(2t_m^2))e^{-it_mH}f\tag 4.13
$$
as $m\to\infty$.

Let constants $\gamma>1$ and $1\ge \theta_1>\rho_j>\theta_j>\rho_N>0$ be fixed such that
$$\align
&\theta_{j-1}\ge \theta_j+\rho_j\q(j=2,\cdots,N-1),\tag 4.14\\
&\max\left\{\gamma(1+\gamma),\frac{2\gamma'_1\gamma'_2}{2-\gamma'_1}\right\}<r_0=\min_{2\le j\le N-1}\{\rho_j/\theta_j\}\tag 4.15
\endalign
$$
for $\gamma'_j$ $(j=1,2)$ defined by (3.9).
Set
$$
\la''_0=\la'_0\theta_{N-1}>0\tag 4.16
$$
with $\la'_0=\la_0-2\kappa$ defined above.
Let $\tau_0>0$ satisfy
$$
0<16\tau_0<\la''_0(<\la'_0<\la_0).\tag 4.17
$$
We take a finite subset $\{{\tilde \la}_\ell\}_{\ell=1}^L$ of $\TT$ such that
$$
\TT\subset \bigcup_{\ell=1}^L({\tilde \la}_\ell-\tau_0,{\tilde\la}_\ell+\tau_0).\tag 4.18
$$
Then we can choose real numbers $\la_\ell\in R^1$, $\tau_\ell>0$ $(\ell=1,2,\cdots,L)$ and $\sig_0>0$ such that
$$\align
&\tau_\ell<\tau_0,\q \sig_0<\tau_0,\q |\la_\ell-{\tilde\la}_\ell|<\tau_0,\\
&\TT\subset\bigcup_{\ell=1}^L(\la_\ell-\tau_\ell,\la_\ell+\tau_\ell),\q (\la_\ell-\tau_\ell,\la_\ell+\tau_\ell)\subset ({\tilde \la}_\ell-\tau_0,{\tilde\la}_\ell+\tau_0),\\
&\text{dist}\{(\la_\ell-\tau_\ell,\la_\ell+\tau_\ell),(\la_k-\tau_k,\la_k+\tau_k)\}>4\sig_0(>0)\q\text{for any}\ \ell\ne k.\tag 4.19
\endalign
$$ 
We note that for $\ell=1,\cdots,L$
$$
\left\{\Lambda \left|\ \tau_\ell\le|\Lambda-(E-\la_\ell)|\le\tau_\ell+4\sig_0\right.\right\}
\cap \Sigma(E)=\emptyset.\tag 4.20
$$
Now let the intervals $\Delta$ and ${\tilde \Delta}$ be so small that
$$
\text{diam}\hskip2pt{\tilde\Delta}<\text{diam}\hskip2pt{\Delta}<{\tilde \tau}_0:=\min_{1\le\ell\le L}\{\sig_0,\tau_\ell\}.\tag 4.21
$$

Returning to (4.6), we have
$$
T_d+E_j\in \text{supp}\hskip2pt\psi,\tag 4.22
$$
if $\psi(T_d+E_j)\ne 0$ in (4.6). By supp\hskip2pt$\psi\subset{\Delta}$, diam\hskip2pt${\Delta}<{\tilde \tau}_0$, and $E\in{\Delta}$, we have from (4.22)
$$
-{\tilde \tau}_0\le T_d-(E-E_j)\le{\tilde\tau}_0.\tag 4.23
$$
Thus we have asymptotically on each state in (4.6)
$$
-2{\tilde\tau}_0\le\frac{|x|^2}{t_m^2}-2(E-E_j)\le2\tilde\tau_0.\tag 4.24
$$
By (4.19), $E_j\in\TT$ is included in just one set $(\la_\ell-\tau_\ell,\la_\ell+\tau_\ell)$ for some $\ell=\ell(j)$ with $1\le \ell(j)\le L$. Since $|E_j-\la_\ellj|<\tau_\ellj$, we have using (4.21)
$$
-2\tau_\ellj-2\sig_0\le\frac{|x|^2}{t_m^2}-2(E-\la_\ellj)\le 2\tau_\ellj+2\sig_0\tag 4.25
$$
on each state in (4.6). Thus 
$$
\sum_{\ell=1}^L \phi_{\sig_0}^2(\left||x|^2/t_m^2-2(E-\la_\ell)\right|<2\tau_\ell+2\sig_0)=1
$$
asymptotically as $m\to\infty$ on (4.6).
Now by the same reasoning that led us to (4.13), we see that (4.3) asymptotically equals as $m\to\infty$
$$
\sum_{\ell=1}^L \sum_{2\le|d|\le N}\phi_{\sig_0}^2(\left||x|^2/t_m^2-2(E-\la_\ell)\right|<2\tau_\ell+2\sig_0)\Psi^2(|x|^2/(2t_m^2))\wtP_d^\Mdm e^{-it_mH}f.\tag 4.26
$$
Since $\phi_{\sig_0}(\left||x|^2/t_m^2-2(E-\la_\ell)\right|<2\tau_\ell+4\sig_0)=1$ on supp $\phi_{\sig_0}(\left||x|^2/t_m^2-2(E-\la_\ell)\right|<2\tau_\ell+2\sig_0)$, (4.26) equals
$$\align
\sum_{\ell=1}^L \sum_{2\le|d|\le N}&\phi_{\sig_0}^2(\left||x|^2/t_m^2-2(E-\la_\ell)\right|<2\tau_\ell+4\sig_0)\\
&\times\phi_{\sig_0}^2(\left||x|^2/t_m^2-2(E-\la_\ell)\right|<2\tau_\ell+2\sig_0)\Psi^2(|x|^2/(2t_m^2))\wtP_d^\Mdm e^{-it_mH}f.\tag 4.27
\endalign
$$
Set
$$
B=\lan x\ran^{-1/2}A\lan x\ran^{-1/2},\q A=\frac{1}{2}(x\cdot p+p\cdot x)=\frac{1}{2}(\lan x,v\ran+\lan v,x\ran).\tag 4.28
$$
We note by Theorem 1.3-(1.23) and (1.25) that on the state $\wtP_d^\Mdm e^{-it_mH}f$
$$
B\sim \sqrt{2T_d}\sim \frac{|x|}{t_m}\tag 4.29
$$
 asymptotically as $t_m\to\infty$. Using this, we replace $\phi_{\sig_0}(\left||x|^2/t_m^2-2(E-\la_\ell)\right|<2\tau_\ell+2\sig_0)$ by $\phi_{\sig_0}(\left|B^2-2(E-\la_\ell)\right|<2\tau_\ell+2\sig_0)$ in (4.27). Let $\varphi(\la)\in C^\infty_0((\sqrt{2(\la_0-2\kappa)}$, $\sqrt{2(\Lambda_0+2\kappa)}))$,
$0\le \varphi(\la)\le 1$, and $\varphi(\la)=1$ on $[\sqrt{2(\la_0-\kappa)},\sqrt{2(\La_0+\kappa)}](\supset$ supp $\Psi(\la^2/2)\cap(0,\infty))$. We insert a factor $\varphi^2(B)$ into (4.27) and then remove the factor $\Psi^2(|x|^2/(2t_m^2))$ using (4.13):
$$
\align
\sum_{\ell=1}^L &\sum_{2\le|d|\le N}\phi_{\sig_0}^2(\left||x|^2/t_m^2-2(E-\la_\ell)\right|<2\tau_\ell+4\sig_0)\\
&\times\phi_{\sig_0}^2(\left|B^2-2(E-\la_\ell)\right|<2\tau_\ell+2\sig_0)
\varphi^2(B)
\wtP_d^\Mdm e^{-it_mH}f.\tag 4.30
\endalign
$$
On supp $\phi_{\sig_0}(\left||x|^2/t_m^2-2(E-\la_\ell)\right|<2\tau_\ell+4\sig_0)$ we have
$$
0<2(E-\la_\ell)-7\tau_0\le\frac{|x|^2}{t_m^2}\le 2(E-\la_\ell)+7\tau_0.\tag 4.31
$$
Since (4.8), $|\la_\ell-\tilde\la_\ell|<\tau_0$ and (4.17) imply
$$
\frac{2(E-\la_\ell)+7\tau_0}{2(E-\la_\ell)-7\tau_0}-1
=\frac{14\tau_0}{2(E-\la_\ell)-7\tau_0}<\frac{14\la''_0/16}{30\la_0/16-\la''_0}
<\theta_{N-1},\tag 4.32
$$
we can apply the partition of unity in Proposition 3.2 to the ring defined by (4.31). Then we obtain
$$\align
e^{-it_mH}f \sim 
&\sum_{\ell=1}^L \sum_{2\le|b|\le N}\sum_{2\le|d|\le N}J_b(x_b/(r_\ell t_m))\phi_{\sig_0}^2(\left||x|^2/t_m^2-2(E-\la_\ell)\right|<2\tau_\ell+4\sig_0)\\
&\times\phi_{\sig_0}^2(\left|B^2-2(E-\la_\ell)\right|<2\tau_\ell+2\sig_0)\varphi^2(B)\wtP_d^\Mdm e^{-it_mH}f,\tag 4.33
\endalign
$$
where 
$$
r_\ell=\sqrt{2(E-\la_\ell)-7\tau_0}>0\q(\ell=1,\cdots,L).\tag 4.34
$$
By the property (3.25), only the terms with $d\le b$ remain in (4.33):
$$\align
e^{-it_mH}f \sim 
&\sum_{\ell=1}^L \sum_{2\le|b|\le N}\sum_{d\le b}J_b(x_b/(r_\ell t_m))\phi_{\sig_0}^2(\left||x|^2/t_m^2-2(E-\la_\ell)\right|<2\tau_\ell+4\sig_0)\\
&\times\phi_{\sig_0}^2(\left|B^2-2(E-\la_\ell)\right|<2\tau_\ell+2\sig_0)\varphi^2(B)\wtP_d^\Mdm e^{-it_mH}f.\tag 4.35
\endalign
$$
Using Theorem 1.3-(1.25), we replace $x_b/t_m$ by $v_b$, and at the same time we introduce a pseudodifferential operator into (4.35):
$$
P_b(t)=\phi_{\sig}(|x_b/t-v_b|^2<u)\tag 4.36
$$
with $u>0$ sufficiently small.
Then (4.35) becomes
$$
\align
e^{-it_mH}f \sim 
&\sum_{\ell=1}^L \sum_{2\le|b|\le N}\sum_{d\le b}P_b^2(t_m)J_b(v_b/r_\ell)\phi_{\sig_0}^2(\left||x|^2/t_m^2-2(E-\la_\ell)\right|<2\tau_\ell+4\sig_0)\\
&\times\phi_{\sig_0}^2(\left|B^2-2(E-\la_\ell)\right|<2\tau_\ell+2\sig_0)\varphi^2(B)\wtP_d^\Mdm e^{-it_mH}f.\tag 4.37
\endalign
$$
We rearrange the order of the factors on the RHS of (4.37) using that the factors mutually commute asymptotically as $m\to\infty$ by Theorem 1.3. Setting
$$\align
G_{b,\la_\ell}(t)=&P_b(t)\phi_{\sig_0}(\left||x|^2/t^2-2(E-\la_\ell)\right|<2\tau_\ell+4\sig_0)\\
&\times\phi_{\sig_0}(\left|B^2-2(E-\la_\ell)\right|<2\tau_\ell+2\sig_0)\varphi(B),\tag 4.38
\endalign
$$
we obtain
$$
e^{-it_mH}f \sim 
\sum_{\ell=1}^L \sum_{2\le|b|\le N}\sum_{d\le b}G_{b,\la_\ell}(t_m)^*
J_b(v_b/r_\ell)G_{b,\la_\ell}(t_m)\wtP_d^\Mdm e^{-it_mH}f.\tag 4.39
$$
Now by some calculus of pseudodifferential operators and Theorem 1.3 we note that $P_b(t)J_b(v_b/r_\ell)$ yields a partition of unity ${\tilde J}_b(x_b/(r_\ell t))$ asymptotically as $m\to\infty$ whose support is close to that of $J_b(x_b/(r_\ell t))$. Then we can recover the terms with $d\not\le b$, and using (1.21), we remove the sum of $\wtP_d^\Mdm$ over $2\le|d|\le N$:
$$
e^{-it_mH}f \sim 
\sum_{\ell=1}^L \sum_{2\le|b|\le N}
G_{b,\la_\ell}(t_m)^*J_b(v_b/r_\ell)G_{b,\la_\ell}(t_m)
e^{-it_mH}f.\tag 4.40
$$
We note that on the RHS, the support with respect to $B^2/2$ of the derivative \linebreak $\phi'_{\sig_0}(\left|B^2-2(E-\la_\ell)\right|<2\tau_\ell+2\sig_0)$ is disjoint with $\Sigma(E)$ by (3.16) and (4.20), and the support of $\varphi'(B)$ is similar by (4.8) and the definition of $\varphi$ above.

We prove the existence of the limit
$$
f_{b,\ell}:=\lim_{t\to\infty}e^{itH}
G_{b,\la_\ell}(t)^*J_b(v_b/r_\ell)G_{b,\la_\ell}(t)
e^{-itH}f\tag 4.41
$$
for $\ell=1,\cdots,L$ and $b$ with $2\le|b|\le N$.

For this purpose we differentiate the function
$$
(e^{itH}
G_{b,\la_\ell}(t)^*J_b(v_b/r_\ell)G_{b,\la_\ell}(t)
e^{-itH}f,g)\tag 4.42
$$
 with respect to $t$, where $f,g\in E_H(\Delta)\HH$. Then writing
$$
D_t^b g(t)=i[H_b,g(t)]+\frac{dg}{dt}(t)\tag 4.43
$$
for an operator-valued function $g(t)$, we have
$$\align
\frac{d}{dt}&
(e^{itH}
G_{b,\la_\ell}(t)^*J_b(v_b/r_\ell)G_{b,\la_\ell}(t)
e^{-itH}f,g)\\
=&(e^{itH}D_t^b(\varphi(B))\phi_{\sig_0}(\left|B^2-2(E-\la_\ell)\right|<2\tau_\ell+2\sig_0)\\
&\q\times \phi_{\sig_0}(\left||x|^2/t^2-2(E-\la_\ell)\right|<2\tau_\ell+4\sig_0)P_b(t)J_b(v_b/r_\ell)G_{b,\la_\ell}(t)e^{-itH}f,g)\\
&+
(e^{itH}\varphi(B)D_t^b\left(\phi_{\sig_0}(\left|B^2-2(E-\la_\ell)\right|<2\tau_\ell+2\sig_0)\right)\\
&\q\times \phi_{\sig_0}(\left||x|^2/t^2-2(E-\la_\ell)\right|<2\tau_\ell+4\sig_0)P_b(t)J_b(v_b/r_\ell)G_{b,\la_\ell}(t)e^{-itH}f,g)\\
&+
(e^{itH}\varphi(B)\phi_{\sig_0}(\left|B^2-2(E-\la_\ell)\right|<2\tau_\ell+2\sig_0)\\
&\q\times D_t^b\left(\phi_{\sig_0}(\left||x|^2/t^2-2(E-\la_\ell)\right|<2\tau_\ell+4\sig_0)\right)P_b(t)J_b(v_b/r_\ell)G_{b,\la_\ell}(t)e^{-itH}f,g)\\
&+
(e^{itH}\varphi(B)\phi_{\sig_0}(\left|B^2-2(E-\la_\ell)\right|<2\tau_\ell+2\sig_0)\\
&\q\times \phi_{\sig_0}(\left||x|^2/t^2-2(E-\la_\ell)\right|<2\tau_\ell+4\sig_0)D_t^b\left(P_b(t)\right)J_b(v_b/r_\ell)G_{b,\la_\ell}(t)e^{-itH}f,g)\\
&+((h.c.)f,g)\\
&+
(e^{itH}i[I_b,G_{b,\la_\ell}(t)
J_b(v_b/r_\ell)G_{b,\la_\ell}(t)]e^{-itH}f,g),\tag 4.44
\endalign
$$
where $(h.c.)$ denotes the adjoint of the operator in the terms preceding it.
\BP

We need the following lemmas (see [K1, Lemmas 4.1 and 4.2]):

\BP

\proclaim{Lemma 4.2} Let Assumption 1.1 be satisfied. Let $E\in R^1-\TT$. Let $F(s)\in C_0^\infty(R^1)$ satisfy $0\le F\le 1$ and the condition that the support with respect to $s^2/2$ of $F(s)$ is disjoint with $\Sigma(E)$. Then there is a constant $d(E)>0$ such that for any interval $\Delta$ around $E$ with {\rm diam} $\Delta<d(E)$, one has
$$
\int_{-\infty}^\infty \left\Vert \frac{1}{\sqrt{\lan x\ran}}F(B)e^{-itH}E_H(\Delta)f\right\Vert^2 dt\le C\Vert f\Vert^2\tag 4.45
$$
for some constant $C>0$ independent of $f\in \HH$.
\endproclaim

\BP

\proclaim{Lemma 4.3} For the pseudodifferential operator $P_b(t)$ defined by (4.36) with $u>0$, there exist norm continuous bounded operators $S(t)$ and $R(t)$ such that
$$
D_t^bP_b(t)=\frac{1}{t}S(t)+R(t)\tag 4.46
$$
and
$$
S(t)\ge 0,\q\Vert R(t)\Vert\le C\lan t\ran^{-2}\tag 4.47
$$
for some constant $C>0$ independent of $t\in R^1$.
\endproclaim

We switch to a smaller interval $\Delta$ if necessary in the followings when we apply Lemma 4.2.

For the first term on the RHS of (4.44) we have
$$
D_t^b(\varphi(B))=\varphi'(B)i[H_b,B]+R_1\tag 4.48
$$
with
$$\align
&\Vert (H+i)^{-1}\lan x\ran^{1/2}i[H_b,B]\lan x\ran^{1/2}(H+i)^{-1}\Vert<\infty,\tag 4.49\\
&\Vert (H+i)^{-1}\lan x\ran R_1 \lan x\ran(H+i)^{-1}\Vert<\infty.\tag 4.50
\endalign
$$
(See section 4 of [K1] for a detailed argument yielding the estimates for the remainder terms $R_1$ here and $S_1(t)$, etc. below.)
By the remark after (4.40), the support with respect to $B^2/2$ of $\varphi'(B)$ is disjoint with $\Sigma(E)$. Hence the condition of Lemma 4.2 is satisfied. Thus using (4.49)-(4.50) and rearranging the order of the factors in the first term on the RHS of (4.44) with some integrable errors, we have by Lemma 4.2:
$$
\text{the 1st term}=(e^{itH}B^{(1)}_2(t)^* B^{(1)}_1(t)e^{-itH}f,g)+(e^{itH}S_1(t)e^{-itH}f,g),\tag 4.51
$$
where $B^{(1)}_j(t)$ $(j=1,2)$ and $S_1(t)$ satisfy
$$\align
&\int_{-\infty}^\infty \Vert B^{(1)}_j(t)e^{-itH}f\Vert^2 dt\le C\Vert f\Vert^2,\tag 4.52\\
&\Vert (H+i)^{-1}S_1(t)(H+i)^{-1}\Vert\le Ct^{-2}\tag 4.53
\endalign
$$ for some constant $C>0$ independent of $f\in E_H(\Delta)\HH$ and $t\in R^1$.

Similarly by another remark after (4.40) and Lemma 4.2, we have a similar bound for the second term on the RHS of (4.44):
$$
\text{the 2nd term}=(e^{itH}B^{(2)}_2(t)^* B^{(2)}_1(t)e^{-itH}f,g)+(e^{itH}S_2(t)e^{-itH}f,g),\tag 4.54
$$
where $B^{(2)}_j(t)$ $(j=1,2)$ and $S_2(t)$ satisfy
$$\align
&\int_{-\infty}^\infty \Vert B^{(2)}_j(t)e^{-itH}f\Vert^2 dt\le C\Vert f\Vert^2,\tag 4.55\\
&\Vert (H+i)^{-1}S_2(t)(H+i)^{-1}\Vert\le Ct^{-2}\tag 4.56
\endalign
$$ for some constant $C>0$ independent of $f\in E_H(\Delta)\HH$ and $t\in R^1$.

For the third term on the RHS of (4.44), we have
$$\align
&\varphi(B)\phi_{\sig_0}(\left|B^2-2(E-\la_\ell)\right|<2\tau_\ell+2\sig_0)\\
&\q\times D_t^b\left(\phi_{\sig_0}(\left||x|^2/t^2-2(E-\la_\ell)\right|<2\tau_\ell+4\sig_0)\right)\\
&=\frac{2}{t}\varphi(B)\phi_{\sig_0}(\left|B^2-2(E-\la_\ell)\right|<2\tau_\ell+2\sig_0)\\
&\q\times \left(\frac{A}{t}-\frac{|x|^2}{t^2}\right)
\phi'_{\sig_0}(\left||x|^2/t^2-2(E-\la_\ell)\right|<2\tau_\ell+4\sig_0)
\\
&\q+S_3(t),\tag 4.57
\endalign
$$
where $S_3(t)$ satisfies
$$
\Vert (H+i)^{-1}S_3(t)(H+i)^{-1}\Vert\le C t^{-2},\q t>1.\tag 4.58
$$
On the support of $\phi'_{\sig_0}(\left||x|^2/t^2-2(E-\la_\ell)\right|<2\tau_\ell+4\sig_0)$, we have 
$$
|x|/t\ge \sqrt{2(E-\la_\ell)-2\tau_\ell-5\sig_0}>0
$$
 by (4.17) and (4.19). Thus there is a large $T>1$ such that for $t\ge T$ we have $|x|>1$ and $\lan x\ran=|x|$, and hence
$$\align
2\left(\frac{A}{t}-\frac{|x|^2}{t^2}\right)&
=\frac{\lan x\ran}{t}\left(\frac{x}{\lan x\ran}\cdot D_x-\frac{|x|}{t}\right)+
\left(D_x\cdot\frac{x}{\lan x\ran}-\frac{|x|}{t}\right)\frac{\lan x\ran}{t}
\\
&=2\frac{\lan x\ran}{t}\left(B-\frac{|x|}{t}\right)+tS_4(t)\tag 4.59
\endalign
$$
with
$\Vert S_4(t)\Vert\le Ct^{-2}$ for $t\ge T$.
By (3.16), we have
$$
\text{supp}\hskip2pt\phi'_{\sig_0}(|s|<2\tau_\ell+4\sig_0)\subset I_1\cup I_2
$$
 with
$$
I_1=[-2\tau_\ell-5\sig_0,-2\tau_\ell-4\sig_0],\q
I_2=[2\tau_\ell+4\sig_0,2\tau_\ell+5\sig_0],\tag 4.60
$$
and 
$$\align
&\phi'_{\sig_0}(|s|<2\tau_\ell+4\sig_0)\ge 0\q \text{for}\ s\in I_1,\tag 4.61\\
&\phi'_{\sig_0}(|s|<2\tau_\ell+4\sig_0)\le 0\q \text{for}\ s\in I_2.\tag 4.62
\endalign
$$
Consider the case $|x|^2/t^2-2(E-\la_\ell)\in I_2$. Then 
$$
\frac{|x|^2}{t^2}\in [2(E-\la_\ell)+2\tau_\ell+4\sig_0,
2(E-\la_\ell)+2\tau_\ell+5\sig_0].\tag 4.63
$$
By the factor $\varphi(B)\phi_{\sig_0}(\left|B^2-2(E-\la_\ell)\right|<2\tau_\ell+2\sig_0)$, we have 
$$
B^2\in [2(E-\la_\ell)-2\tau_\ell-3\sig_0,2(E-\la_\ell)+2\tau_\ell+3\sig_0]\tag 4.64
$$
 and
$B\ge \sqrt{2\la'_0}>0$. Thus
$$
B-\frac{|x|}{t}\le 0.\tag 4.65
$$
Therefore
by (4.59) and (4.62), (4.57) is positive in this case up to an integrable error. Similarly we see that (4.57) is positive also in the case
 $|x|^2/t^2-2(E-\la_\ell)\in I_1$. Rearranging the order of the factors in the third term on the RHS of (4.44) with an integrable error, we see that it has the form
$$
\text{the 3rd term}=(e^{itH}A(t)^*A(t)e^{-itH}f,g)+(e^{itH}S_5(t)e^{-itH}f,g)\tag 4.66
$$
with
$$
\Vert (H+i)^{-1}S_5(t)(H+i)^{-1}\Vert\le Ct^{-2}.\tag 4.67
$$

The fourth term on the RHS of (4.44) has a similar form by virtue of Lemma 4.3.

The fifth term $((h.c.)f,f)$ is treated similarly to the terms above.

The sixth term on the RHS of (4.44) satisfies
$$
|\text{the 6th term}|\le Ct^{-1-\min\{\ep,\ep_1\}}\Vert f\Vert \Vert g\Vert.\tag 4.68
$$
This estimate follows if we note with using (3.25) and some calculus of pseudodifferential operators as stated after (4.39) that the factor $G_{b,\la_\ell}(t)^*J_b(v_b/r_\ell)G_{b,\la_\ell}(t)$ restricts the coordinates in the region: $|x_\al|^2>\rho_{|b|}|x|^2/2$.

Summarizing we have proved that (4.44) is written as
$$\align
\frac{d}{dt}&
(e^{itH}
G_{b,\la_\ell}(t)^*J_b(v_b/r_\ell)G_{b,\la_\ell}(t)
e^{-itH}f,g)\\
&=(e^{itH}A(t)^*A(t)e^{-itH}f,g)+\sum_{k=1}^2(e^{itH}B^{(k)}_2(t)^*B^{(k)}_1(t)e^{-itH}f,g)+(S_6(t)f,g),\tag 4.69
\endalign
$$
where with some constant $C>0$ independent of $t>T$ and $f\in \HH$
$$
\align
&\int_T^\infty\Vert B^{(k)}_j(t)e^{-itH}E_H(\Delta)f\Vert^2\le C\Vert f\Vert^2,\q (j,k=1,2)\tag 4.70\\
&\Vert (H+i)^{-1}S_6(t)(H+i)^{-1}\Vert\le Ct^{-1-\min\{\ep,\ep_1\}}.\tag 4.71
\endalign
$$

Integrating (4.69) with respect to $t$ on an interval $[T_1,T_2]\subset[T,\infty)$, we obtain
$$\align
& \left.(e^{itH}
G_{b,\la_\ell}(t)^*J_b(v_b/r_\ell)G_{b,\la_\ell}(t)
e^{-itH}f,g)\right|_{t=T_1}^{T_2}\\
&\q=
\int_{T_1}^{T_2}(A(t)e^{-itH}f,A(t)e^{-itH}g) dt\\
&\qq+
\sum_{k=1}^2\int_{T_1}^{T_2}(B^{(k)}_1(t)e^{-itH}f,B^{(k)}_2(t)e^{-itH}g)dt+\int_{T_1}^{T_2}(S_6(t)f,g)dt.\tag 4.72
\endalign
$$
Hence using (4.70), (4.71) and the uniform boundedness of $G_{b,\la_\ell}(t)$ in $t>1$, we have
$$
\int_{T_1}^{T_2}\Vert A(t)e^{-itH}g\Vert^2 dt\le C\Vert g\Vert^2\tag 4.73
$$
for some constant $C>0$ independent of $T_2>T_1\ge T$ and $g\in E_H(\Delta)\HH$.

\define\tf{{\tilde f}}

(4.73) and (4.72) with (4.70) and (4.71) then yield that
$$
\left|\left.(e^{itH}
G_{b,\la_\ell}(t)^*J_b(v_b/r_\ell)G_{b,\la_\ell}(t)
e^{-itH}f,g)\right|_{t=T_1}^{T_2}\right|
\le \de(T_1)\Vert f\Vert\Vert g\Vert\tag 4.74
$$
for some $\de(T_1)>0$ with $\de(T_1)\to 0$ as $T_2>T_1\to\infty$. This means that
the limit
$$
\tf_{b}^1=\lim_{t\to\infty}\sum_{\ell=1}^L e^{itH}
G_{b,\la_\ell}(t)^*J_b(v_b/r_\ell)G_{b,\la_\ell}(t)
e^{-itH}f\tag 4.75
$$
exists for any $f\in E_H(\Delta)\HH$ and $b$ with $2\le|b|\le N$ if $\Delta$ is an interval sufficiently small around $E\in R^1-\TT$: diam\hskip2pt$\Delta<d(E)$. Then the asymptotic decomposition (4.40) implies
$$
f=\sum_{2\le|b|\le N} \tf_b^1\tag 4.76
$$
for $f=\psi(H)f=E_H(\Delta)f$.
Further by the existence of the limit (4.75) and $f=E_H(\Delta)f$, we see that $\tf_b^1$ satisfies
$$
E_H(\Delta)\tf_b^1=\tf_b^1\tag 4.77
$$
in a way similar to the proof of the intertwining property of wave operators.

\MP

Now returning to the first $\Phi(H)f$, and noting that supp\hskip2pt$\Phi$ is compact in $R^1-\TT$, we take a finite number of open intervals $\Delta_{j_0}\Subset R^1-\TT$ such that $E_{j_0}\in\Delta_{j_0}$, diam\hskip2pt$\Delta_{j_0}<d(E_{j_0})$, and supp\hskip2pt$\Phi\Subset\bigcup_{j_0}^{\text{finite}}\Delta_{j_0}\Subset R^1-\TT$. Then we can take $\psi_{j_0}\in C_0^\infty(\Delta_{j_0})$ such that $\Phi(H)f=\sum_{j_0}^{\text{finite}}\psi_{j_0}(H)f$. Thus from (4.75)-(4.77), we obtain the existence of the limit for $2\le|b|\le N$:
$$
\tf_b^1=\lim_{t\to\infty}\sum_{j_0}^{\text{finite}} \sum_{\ell=1}^L e^{itH}
G_{b,\la_\ell}(t)^*J_b(v_b/r_\ell)G_{b,\la_\ell}(t)
e^{-itH}\psi_{j_0}(H)f,
\tag 4.78
$$
and the relations
$$
\Phi(H)f=\sum_{2\le|b|\le N}\tf_b^1,\q E_H(\Delta)\tf_b^1=\tf_b^1\tag 4.79
$$
for any set $\Delta\Subset R^1-\TT$ with supp\hskip2pt$\Phi\subset\Delta$.

Set
$$
\sig_j=\sqrt{\gamma^{-1}\rho_j\la'_0/2},\q
\de_j=\sqrt{\gamma\theta_j\La'_0}\q (j=2,3,\cdots,N,\q \theta_N=0).\tag 4.80
$$
Then by (4.78), some calculus of pseudodifferential operators, and
$$
\text{supp}\hskip2pt{\left(J_b(x_b/r_\ell)\phi_{\sig_0}(\left||x|^2-2(E-\la_\ell)\right|<2\tau_\ell+4\sig_0)\right)}\Subset {\tilde T}_b(\gamma^{-1}\rho_{|b|},\gamma\theta_{|b|}),\tag 4.81
$$
which follows from (3.19)-(3.21), we see that as $t\to\infty$
$$\align
e^{-itH}\tf_b^1&\sim 
\sum_{k=1}^K \sum_{\ell=1}^L G_{b,\la_\ell}(t)^*J_b(v_b/r_\ell)G_{b,\la_\ell}(t)
 e^{-itH}E_H(\Delta_k)f\\
&\sim
\prod_{\al\not\leq b}F(|x_\al|\ge \sig_{|b|} t)F(|x^b|\le \de_{|b|} t)\\
&\qq\qq\times \sum_{k=1}^K\sum_{\ell=1}^LG_{b,\la_\ell}(t)^*J_b(v_b/r_\ell)G_{b,\la_\ell}(t)
 e^{-itH}E_H(\Delta_k)f\\
&\sim \prod_{\al\not\leq b}F(|x_\al|\ge \sig_{|b|} t)F(|x^b|\le \de_{|b|} t) e^{-itH}\tf_b^1.\tag 4.82
\endalign
$$
(4.79) and (4.82) imply
$$
\tf_b^1\in S_b^{1\sig_{|b|}\de_{|b|}}(\Delta)\q(2\le|b|\le N).\tag 4.83
$$

\BP


\F
II) A refinement:
\BP

\define\deL{{\de'}}
\define\sigL{{\sig'}}

As in (4.18)-(4.19), we take a finite subset $\{{\tilde\la}^b_\ell\}_{\ell=1}^{L_b}$ of $\TT_b$ for a constant $\tau_0^b>0$ with $\tau_0^b<\tau_0$ such that
$$
\TT_b\subset\bigcup_{\ell=1}^{L_b}({\tilde\la}^b_\ell-\tau_0^b,{\tilde\la}^b_\ell+\tau^b_0),\tag 4.84
$$ and choose real numbers $\la_\ell^b\in R^1$, $\tau_\ell^b>0$ $(\ell=1,\cdots,L_b)$ and $\sig_0^b>0$ such that
$$\align
&\tau_\ell^b<\tau_0^b,\q \sig_0^b<\tau_0^b,\q |\la_\ell^b-{\tilde\la}^b_\ell|<\tau_0^b,\\
&\TT_b\subset\bigcup_{\ell=1}^{L_b}(\la^b_\ell-\tau^b_\ell,\la^b_\ell+\tau^b_\ell),\q (\la^b_\ell-\tau^b_\ell,\la^b_\ell+\tau^b_\ell)\subset ({\tilde \la}^b_\ell-\tau_0^b,{\tilde\la}^b_\ell+\tau_0^b),\\
&\text{dist}\{(\la^b_\ell-\tau^b_\ell,\la^b_\ell+\tau^b_\ell),(\la^b_k-\tau^b_k,\la_k^b+\tau_k^b)\}>4\sig_0^b(>0)\q\text{for any}\ \ell\ne k.\tag 4.85
\endalign
$$ 
We set $\TT^{F}_b=\{\la_\ell\}_{\ell=1}^{L_b}$ and
$$
{\tilde\tau}_0^b=\min_{1\le\ell\le L_b}\{\sig_0^b,\tau_\ell^b\}.\tag 4.86
$$
Then, we take $\psi_1(\la)\in C_0^\infty(R^1)$ such that
$$\align
&0\le\psi_1\le 1,\tag 4.87\\
&\psi_1(\la)=
\cases
1&\q  \text{for any}\ \la\ \text{with}\ |\la-\la_\ell|\le {\tilde\tau}_0^b/2\ \text{for some}\ \la_\ell\in\TT^F_b\\
0&\q \text{for any}\ \la\ \text{with}\ |\la-\la_\ell|\ge {\tilde\tau}_0^b\ \text{for all}\ \la_\ell\in\TT^F_b
\endcases
\tag 4.88
\endalign
$$
and
we divide (4.78) as follows:
$$
\tf_b^1=h_b+g_b,\tag 4.89
$$
where
$$\align
h_b&=\lim_{t\to\infty}\sum_{j_0}^{\text{finite}} \sum_{\ell=1}^L e^{itH}
G_{b,\la_\ell}(t)^*\psi_1(H^b)J_b(v_b/r_\ell)G_{b,\la_\ell}(t)
e^{-itH}\psi_{j_0}(H)f,\tag 4.90\\
g_b&=\lim_{t\to\infty}\sum_{j_0}^{\text{finite}} \sum_{\ell=1}^L e^{itH}
G_{b,\la_\ell}(t)^*(I-\psi_1)(H^b)J_b(v_b/r_\ell)G_{b,\la_\ell}(t)
e^{-itH}\psi_{j_0}(H)f.
\tag 4.91
\endalign
$$
The proof of the existence of these limits is similar to that of ${\tilde f}_b^1$ in (4.78), since the change in the present case is the appearance of the commutator $[H,\psi_1(H^b)]=[I_b,\psi(H^b)]$ whose treatment is quite the same as  that of the commutators including $I_b$ in (4.68).
We introduce the decomposition (1.21) into $h_b$ and $g_b$ on the left of $e^{-itH}\psi_{j_0}(H)f$ as in (4.3).
Then by the factor $J_b(v_b/r_\ell)G_{b,\la_\ell}(t)$, we see that only the terms with $d\le b$ in the sum in (4.3) remain asymptotically as $t=t_m\to\infty$ by the arguments similar to step I). On each summand $P_{d,E_j}\whP_{|d|-1}^{{\widehat M}_{d}^m}$ in these terms (see (4.5)), $H^b$ asymptotically equals $H^b_d=T_d^b+H^d=T_d^b+E_j\sim|x^b_d|^2/(2t_m^2)+E_j\sim|x^b|^2/(2t_m^2)+E_j$, where for $d\le b$, $H^b_d=T^b_d+H^d=H_d-T_b$, $T^b_d=T_d-T_b$ and $x^b=(x^b_d,x^d)$ is a clustered Jacobi coordinate inside the coordinate $x^b$. 
 Thus we have
$$\align
&\psi_1(H^b)J_b(v_b/r_\ell)G_{b,\la_\ell}(t)P_{d,E_j}\whP_{|d|-1}^{{\widehat M}_{d}^m}e^{-it_mH}\psi_{j_0}(H)f\\
&\sim\psi_1(|x^b|^2/(2t_m^2)+E_j)J_b(v_b/r_\ell)G_{b,\la_\ell}(t)P_{d,E_j}\whP_{|d|-1}^{{\widehat M}_{d}^m}e^{-it_mH}\psi_{j_0}(H)f\tag 4.92
\endalign
$$
as $m\to\infty$. If $\psi_1(|x^b|^2/(2t_m^2)+E_j)\ne 0$, then for some $\ell=1,\cdots,L_b$
$$
\left|\frac{|x^b|^2}{2t_m^2}-(\la_\ell^b-E_j)\right|\le {\tilde\tau}_0^b.\tag 4.93
$$
If $\ell$ is a (unique) $\ellj$ such that $E_j\in(\la_\ellj^b-\tau_\ellj^b,\la_\ellj+\tau_\ellj^b)$, we have
$$
\frac{|x^b|^2}{2t_m^2}\le\tau_\ellj^b+{\tilde\tau}_0^b<\tau_0^b+{\tilde\tau}_0^b.\tag 4.94
$$
Thus setting $\deL=\sqrt{2(\tau_0^b+{\tilde\tau}_0^b)}$, we have
$$
|x^b|\le\deL t_m.\tag 4.95
$$
If $\ell\ne\ellj$, we have by (4.93)
$$
0\le \la_\ell^b-E_j+{\tilde\tau}_0^b,
$$
from which and (4.85)-(4.86) follows
$$
\la_\ell^b-E_j\ge 4\sig_0^b.
$$
Thus from (4.93)
$$
\frac{|x^b|^2}{2t_m^2}\ge 4\sig_0^b-{\tilde\tau}_0^b\ge 3\sig_0^b\ge 3{\tilde\tau}_0^b.\tag 4.96
$$
Setting $\sigL=\sqrt{6{\tilde\tau}_0^b}$ we then have for $\ell\ne \ellj$
$$
|x^b|\ge \sigL t_m.\tag 4.97
$$
Therefore $h_b$ can be decomposed as
$$
h_b=f_{b}^\deL+g_{b1}^\sigL,\tag 4.98
$$
where
$$
\align
f_{b}^\deL&=\lim_{t\to\infty}\sum_{j_0}^{\text{finite}} \sum_{\ell=1}^L e^{itH}
G_{b,\la_\ell}(t)^*F(|x^b|\le\deL t)\psi_1(H^b)J_b(v_b/r_\ell)G_{b,\la_\ell}(t)
e^{-itH}\psi_{j_0}(H)f,\tag 4.99\\
g_{b1}^\sigL&=\lim_{t\to\infty}\sum_{j_0}^{\text{finite}} \sum_{\ell=1}^L e^{itH}
G_{b,\la_\ell}(t)^*F(|x^b|\ge \sigL t)\psi_1(H^b)J_b(v_b/r_\ell)G_{b,\la_\ell}(t)
e^{-itH}\psi_{j_0}(H)f,\tag 4.100
\endalign
$$
The existence of the limit (4.99) is proved similarly to that of (4.90) by rewriting the factor $F(|x^b|\le\deL t)$ as a smooth one and absorbing it into $J_b(v_b/r_\ell)G_{b,\la_\ell}(t)$ with changing the constants in it suitably. The existence of (4.100) then follows from this, (4.78) and (4.98).

For $g_b$, similarly to $g_{b1}^\sigL$ we obtain
$$
g_b
=\lim_{t\to\infty}\sum_{j_0}^{\text{finite}} \sum_{\ell=1}^L e^{itH}
G_{b,\la_\ell}(t)^*F(|x^b|\ge\sigL t)(I-\psi_1)(H^b)J_b(v_b/r_\ell)G_{b,\la_\ell}(t)
e^{-itH}\psi_{j_0}(H)f.\tag 4.101
$$
Setting
$$\align
g_b^\sigL&=g_{b1}^\sigL+g_b\\
&=\lim_{t\to\infty}\sum_{j_0}^{\text{finite}} \sum_{\ell=1}^L e^{itH}
G_{b,\la_\ell}(t)^*F(|x^b|\ge\sigL t)J_b(v_b/r_\ell)G_{b,\la_\ell}(t)
e^{-itH}\psi_{j_0}(H)f,\tag 4.102
\endalign
$$
we obtain a decomposition of ${\tilde f}_b^1$:
$$
{\tilde f}_b^1=f_b^\deL+g_b^\sigL,\tag 4.103
$$
where $f_b^\deL$ and $g_b^\sigL$ satisfy
$$\align
e^{-itH}f_b^\deL&\sim 
\prod_{\al\not\leq b}F(|x_\al|\ge \sig_{|b|} t)F(|x^b|\le \deL t)e^{-itH}f_b^\deL,\tag 4.104\\
e^{-itH}g_b^\sigL&\sim 
\prod_{\al\not\leq b}F(|x_\al|\ge \sig_{|b|} t)F(|x^b|\le \de_{|b|} t)F(|x^b|\ge\sigL t)e^{-itH}g_b^\sigL.\tag 4.105
\endalign
$$
We can prove the existence of the limits
$$
f_b^1=\lim_{\deL\downarrow0}f_b^\deL,\q
g_b^1=\lim_{\sigL\downarrow0}g_b^\sigL\tag 4.106
$$
in the same way as in Enss [E3, Lemma 4.8], because we can take $\psi_1$ in (4.87)-(4.88) monotonically decreasing when ${\tilde\tau}_0^b\downarrow0$ and the factors $F(|x^b|\le \deL t)$ and $F(|x^b|\ge\sigL t)$ can be treated similarly to $\psi_1$ by regarding $x^b/t$ as a single variable. Further we have as in (4.77)
$$
E_H(\Delta)f_b^1=f_b^1,\q E_H(\Delta)g_b^1=g_b^1,\tag 4.107
$$
which, (4.104) and (4.106) imply
$$
f_b^1 \in S_b^1.\tag 4.108
$$

Thus we have a decomposition:
$$
{\tilde f}_b^1=f_b^1+g_b^1,\q f_b^1\in S_b^1.\tag 4.109
$$
$g_b^\sigL$ can be decomposed further by using the partition of unity of the ring $\sigL\le |x^b|/t\le \de_{|b|}$ with regarding $x^b$ as a total variable $x$ in Proposition 3.2. Arguing similarly to steps I) and II), we can prove that $g_b^1$ can be decomposed as a sum of the elements $f_d^1$ of $S_d^1$ with $d< b$.
Combining this with (4.76), (4.108) and (4.109), we obtain (4.1).
\ $\square$
\enddemo
\BP

We remark that Theorem 4.1 implies the asymptotic completeness when the long-range part $V_\al^L$ vanishes for all pairs $\al$, because in this case we see straightforwardly that $S_b^1=\RR(W_b^\pm)$, where $W_b^\pm$ are the short-range wave operators defined by
$$
W_b^\pm=\text{s-}\lim_{t\to\pm\infty}e^{itH}e^{-itH_b}P_b.\tag 4.110
$$
For the case when long-range part does not vanish, 
we have the following 

\proclaim{Theorem 4.2} Let Assumptions 1.1 and 1.2 be satisfied. Then 
\MP

\F
i) For $2(2+\ep)^{-1}<r\le1$
$$
S_b^r=S_b^1.\tag 4.111
$$
\MP

\F
ii) If $\ep>2(2+\ep)^{-1}$, i.e. when $\ep>\sqrt{3}-1$, we have for all $r$ with $0\le r\le 1$
$$
S_b^r=S_b^1.\tag 4.112
$$
\MP

\F
iii) If $\ep>1/2$ and $V^L_\al(x_\al)\ge 0$ for all pairs $\al$, then we have for all $r$ with $0\le r\le 1$
$$
S_b^r=S_b^1.\tag 4.113
$$
\endproclaim
\demo{Proof} i) and ii) follow from Proposition 5.8 of [D] and Proposition 2.3 above. (4.112) for $r=0$ follows from the proof of Proposition 5.8 of [D]. iii) follows from Theorem 1.1 and Proposition 4.3 of [K2] and (4.116) below: Note that $\RR(\Omega_b^\psi)$ in Theorem 1.1-(1.31) of [K2] constitutes a dense subset of $S_b^1$, when $\psi$ varies in $C_0^\infty(R^1-\TT)$.
\ $\square$
\enddemo

From Theorem 4.2-ii), iii) and Theorem 4.1 follows

\proclaim{Theorem 4.3} Let Assumptions 1.1 and 1.2 be satisfied with
 $\ep>2(2+\ep)^{-1}$ or with $\ep>1/2$ and $V^L_\al(x_\al)\ge 0$ for all pairs $\al$. Then we have for all $r$ with $0\le r\le 1$
$$
\bigoplus_{2\le|b|\le N} S_b^r =\HH_c(H).\tag 4.114
$$
\endproclaim

In the next section, we will construct modified wave operators: 
$$
W_b^\pm=\text{s-}\lim_{t\to\pm\infty}e^{itH}J_be^{-itH_b}P_b\tag 4.115
$$
with $J_b$ being an extension of $J$ of [IK] to the $N$-body case.
We will then prove
$$
\RR(W_b^\pm)=S_b^0,\tag 4.116
$$
which and Theorems 4.2 and 4.3 imply

\proclaim{Theorem 4.4} Let Assumptions 1.1 and 1.2 be satisfied with
 $\ep>2(2+\ep)^{-1}$ or with $\ep>1/2$ and $V^L_\al(x_\al)\ge 0$ for all pairs $\al$. Then we have for all $r$ with $0\le r\le 1$
$$
\RR(W_b^\pm)=S_b^r,\tag 4.117
$$
and
$$
\bigoplus_{2\le|b|\le N}\RR(W_b^\pm)=\HH_c(H).\tag 4.118
$$
\endproclaim

One might expect that (4.114) and (4.118) are always true, but it is denied:

\proclaim{Theorem 4.5} Let Assumptions 1.1 and 1.2 be satisfied and let $N\ge 3$. 
Then the followings hold:
\MP

\F
i) Let $2\le|b|\le N$ and
let $E_b(r)$ be the orthogonal projection onto $S_b^r$ $(0\le r\le 1)$. Then $E_b(r_1)\le E_b(r_2)$ for $0\le r_1\le r_2\le 1$, and the discontinuous points of $E_b(r)$ with respect to $r\in[0,1]$ in the strong operator topology are at most countable.
\MP

\F
ii) Let $0<\ep<1/2$ in Assumption 1.1. Then there are long-range pair potentials $V_\al(x_\al)$ such that for some cluster decomposition $b$ with $2\le|b|\le N$, $E_b(r)$ is discontinuous at $r=r_0$, where $\ep<r_0:=(\ep+1)/3<1/2$. In particular, there are real numbers $r_1$ and $r_2$ with $0\le r_1<r_0<r_2\le 1$ such that 
$$
S_b^{r_1}\subsetneqq S_b^{r_2}.\tag 4.119
$$
\endproclaim
\demo{Proof} i) By Proposition 2.3, $S_b^r$ $(0\le r\le 1)$ is a family of closed subspaces of a separable Hilbert space $\HH$ that increases when $r\in[0,1]$ increases. Thus the corresponding orthogonal projection $E_b(r)$ $(0\le r\le 1)$ onto $S_b^r$ increases as $r$ increases, and hence has at most a countable number of discontinuous points with respect to $r\in[0,1]$ in the strong operator topology.

ii) holds by Theorem 4.3 of [Y], Theorem 4.1 and Proposition 2.3, for $b$, $|b|=N-1$, with a suitable choice of pair potentials that satisfy Assumption 1.1.  In fact, the sum of the ranges $\RR(W_n)$ of Yafaev's wave operators $W_n$ $(n=1,2,\cdots)$ in Theorem 4.3 of [Y] constitutes a subspace of $(E_b(r_0+0)-E_b(r_0-0))\HH$ for $b$ with $|b|= N-1$ by his construction of $W_n$, which means that $E_b(r)$ is discontinuous at $r=r_0$. Here $E_b(r_0\pm 0)=\text{s-}\lim_{r\to r_0\pm 0}E_b(r)$. \ $\square$
\enddemo

\vskip 12pt

\BP

\subhead\nofrills
5. A Characterization of the Ranges of Wave Operators
\endsubhead
\par
\vskip 12pt

The purpose in this section is to prove relation (4.116) for general long-range pair potentials $V_\alpha(x_\alpha)$ under Assumptions 1.1 and 1.2. The inclusion
$$
\RR(W_b^\pm)\subset S_b^0\tag 5.1
$$
is a trivial relation for any form of definition of the wave operators $W_b^\pm$. Thus our main concern is to prove the reverse inclusion
$$
S_b^0\subset \RR(W_b^\pm).\tag 5.2
$$
The proof of this inclusion is essentially the same for any definition of wave operators and is not difficult in the light of Enss method [E1]. As announced, we here consider the wave operators of the form
$$
W_b^\pm=\text{s-}\lim_{t\to\pm\infty}e^{itH}J_be^{-itH_b}P_b,\tag 5.3
$$
where $J_b$ is an extension of the identification operator or stationary modifier introduced in [IK] for two-body long-range case. The first task in this section is to construct $J_b$. Our relation (5.2) then follows from the definition of the scattering spaces $S_b^0$ and properties of $J_b$ by Enss method.

To make the descriptions simple we hereafter consider the case $V_\alpha^S=0$ for all pairs $\alpha$. The recovery of the short-range potentials in the following arguments is easy.

Let a $C^\infty$ function $\chi_0(x)$ of $x\in R^\nu$ satisfy
$$
\chi_0(x)=\cases
1 & (|x|\ge 2)\\
0 & (|x|\le 1).
\endcases \tag 5.4
$$
To define $J_b$ we introduce time-dependent potentials $I_{b\rho}(x_b,t)$ for $\rho\in(0,1)$:
$$
I_{b\rho}(x_b,t)=I_b(x_b,0)\prod_{k=1}^{k_b}\chi_0(\rho z_{bk})\chi_0(\langle \log\langle t\rangle\rangle z_{bk}/\langle t\rangle).\tag 5.5
$$
Then $I_{b\rho}(x_b,t)$ satisfies
$$
|\partial_{x_b}^\beta I_{b\rho}(x_b,t)|\le C_\beta \rho^{\epsilon_0}\langle t\rangle^{-\ell}\tag 5.6
$$
for any $\ell\ge 0$ and $0<\epsilon_0<\epsilon$ with $\epsilon_0+\ell<|\beta|
+\epsilon$, where $C_\beta>0$ is a constant independent of $t,x_b$ and $\rho$.

Then we can apply almost the same arguments as in section 2 of [IK] to get a solution $\varphi_b(x_b,\xi_b)$ of the eikonal equation:
$$
\frac{1}{2}|\nabla_{x_b}\varphi_b(x_b,\xi_b)|^2+I_b(x_b,0)=\frac{1}{2}|\xi_b|^2
\tag 5.7
$$
in some conic region in phase space. More exactly we have the following theorems. Let
$$
\cos(z_{bk},\zeta_{bk}):=\frac{z_{bk}\cdot \zeta_{bk}}{|z_{bk}|_e|\zeta_{bk}|_e},
$$
where $|z_{bk}|_e= (z_{bk}\cdot z_{bk})^{1/2}$ is the Euclidean norm. We then set for $R_0,d>0$ and $\theta\in(0,1)$
$$
\Gamma_{\pm}(R_0,d,\theta)=\{(x_b,\xi_b)\ |\ |z_{kb}|\ge R_0,|\zeta_{bk}|\ge d, \pm\cos(z_{bk},\zeta_{bk})\ge \theta\ (k=1,\cdots,k_b)\},
$$
where $\zeta_{bk}$ is the variable conjugate to $z_{bk}$.

\proclaim{Theorem 5.1} Let Assumption 1.1 be satisfied with $V_\alpha^S=0$ for all pairs $\alpha$. Then there exists a $C^\infty$ function $\phi_b^\pm(x_b,\xi_b)$ that satisfies the following properties: For any $0<\theta,d<1$, there exists a constant $R_0> 1$ such that for any $(x_b,\xi_b)\in\Gamma_{\pm}(R_0,d,\theta)$
$$
\frac{1}{2}|\nabla_{x_b}\phi_b^\pm(x_b,\xi_b)|^2+I_b(x_b,0)=\frac{1}{2}|\xi_b|^2
\tag 5.8
$$
and
$$
|\partial_{x_b}^\alpha\partial_{\xi_b}^\beta(\phi_b^\pm(x_b,\xi_b)- x_b\cdot\xi_b)|\le
\cases
 C_{\alpha\beta}\left(\max_{1\le k\le k_b}\langle z_{bk}\rangle\right)^{1-\epsilon},& \alpha=0\\
 C_{\alpha\beta}\left(\min_{1\le k\le k_b}\langle z_{bk}\rangle\right)^{1-\epsilon-|\alpha|},& \alpha\ne 0,
\endcases
\tag 5.9
$$
where $C_{\alpha\beta}>0$ is a constant independent of $(x_b,\xi_b)\in\Gamma_\pm(R_0,d,\theta)$.
\endproclaim

From this we can derive the following theorem in quite the same way as that for Theorem 2.5 of [IK]. Let $0<\theta<1$ and let $\psi_\pm(\tau)\in C^\infty([-1,1])$ satisfy
$$
\split
&0\le \psi_\pm(\tau)\le 1,\\
&\psi_+(\tau)=
\cases 1 &\text{for}\ \theta\le \tau\le 1,\\
0&\text{for}\ -1\le\tau\le\theta/2,
\endcases\\
&\psi_-(\tau)=
\cases 0 &\text{for}\ -\theta/2\le\tau\le 1,\\
1&\text{for}\ -1\le \tau\le-\theta.
\endcases
\endsplit
$$
We set
$$
\chi_\pm(x_b,\xi_b)=\prod_{k=1}^{k_b}\psi_\pm(\cos(z_{bk},\zeta_{bk}))
$$
and define $\varphi_b(x_b,\xi_b)=\varphi_{b,\theta,d,R_0}(x_b,\xi_b)$ by
$$
\split
\varphi_b(x_b,\xi_b)=&\{(\phi_b^+(x_b,\xi_b)-x_b\cdot\xi_b)\chi_+(x_b,\xi_b)+
(\phi_b^-(x_b,\xi_b)-x_b\cdot\xi_b)\chi_-(x_b,\xi_b)\}\\
&\times\prod_{k=1}^{k_b}\chi_0(2\zeta_{bk}/d)\chi_0(2z_{bk}/R_0)+x_b\cdot\xi_b
\endsplit
$$
for
$d,R_0>0$. Note that $\varphi_{b,\theta,d,R_0}(x_b,\xi_b)=\varphi_{b,\theta,d',R_0'}(x_b,\xi_b)$ when $|z_{bk}|\ge\max(R_0,R_0')$, $|\zeta_{bk}|\ge\max(d,d')$ for all $k$. We then have

\proclaim{Theorem 5.2} Let Assumption 1.1 be satisfied with $V_\alpha^S=0$ for all pairs $\alpha$. Let $0<\theta<1$ and $d>0$. Then there exists a constant $R_0>1$ such that the $C^\infty$ function $\varphi_b(x_b,\xi_b)$ defined above satisfies the following properties.
\SP

\F
i) For $(x_b,\xi_b)\in \Gamma_+(R_0,d,\theta)\cup\Gamma_-(R_0,d,\theta)$, $\varphi_b$ is a solution of
$$
\frac{1}{2}|\nabla_{x_b}\varphi_b(x_b,\xi_b)|^2+I_b(x_b,0)=\frac{1}{2}|\xi_b|^2.
\tag 5.10
$$

\F
ii) For any $(x_b,\xi_b)\in R^{2\nu(|b|-1)}$ and multi-indices $\alpha,\beta$, $\varphi_b$ satisfies
$$
|\partial_{x_b}^\alpha\partial_{\xi_b}^\beta(\varphi_b(x_b,\xi_b)-x_b\cdot\xi_b)|
\le
\cases
 C_{\alpha\beta}\left(\max\langle z_{bk}\rangle\right)^{1-\epsilon}, &\alpha=0\\
 C_{\alpha\beta}\left(\min\langle z_{bk}\rangle\right)^{1-\epsilon-|\alpha|}, &\alpha\ne0.
\endcases\tag 5.11
$$
In particular, if $\alpha\ne 0$,
$$
|\partial_{x_b}^\alpha\partial_{\xi_b}^\beta(\varphi_b(x_b,\xi_b)-x_b\cdot\xi_b
)|\le C_{\alpha\beta}R_0^{-\epsilon_0}
\left(\min\langle z_{bk}\rangle\right)^{1-\epsilon_1-|\alpha|}\tag 5.12
$$
for any $\epsilon_0,\epsilon_1\ge 0$ with $\epsilon_0+\epsilon_1=\epsilon$. Further
$$
\varphi_b(x_b,\xi_b)=x_b\cdot\xi_b\q \text{when } |z_{bk}|\le R_0/2\  \text{or}\ |\zeta_{bk}|\le d/2\ \text{for some}\ k.\tag 5.13
$$

\F
iii) Let
$$
a_b(x_b,\xi_b)=e^{-i\varphi_b(x_b,\xi_b)}\left(T_b+I_b(x_b,0)-\frac{1}{2}|\xi_b|^2\right) e^{i\varphi_b(x_b,\xi_b)}.\tag 5.14
$$
Then
$$
a_b(x_b,\xi_b)=\frac{1}{2}|\nabla_{x_b}\varphi_b(x_b,\xi_b)|^2+I_b(x_b,0)-\frac{1}{2}|\xi_b|^2+i(T_b\varphi_b)(x_b,\xi_b)\tag 5.15
$$
and
$$
|\partial_{x_b}^\alpha\partial_{\xi_b}^\beta a_b(x_b,\xi_b)|
\le
\cases
C_{\alpha\beta}\left(\min\langle z_{bk}\rangle\right)^{-1-\epsilon-|\alpha|},& (x_b,\xi_b)\in\Gamma_+(R_0,d,\theta)\cup\Gamma_-(R_0,d,\theta)\\
C_{\alpha\beta}\left(\min\langle z_{bk}\rangle\right)^{-\epsilon-|\alpha|}\langle\xi_b\rangle,&
\text{otherwise}.
\endcases\tag 5.16
$$
\endproclaim

We now define $J_b=J_{b,\theta,d,R_0}$ by
$$
J_bf(x_b)=(2\pi)^{-\nu(|b|-1)}\int_{R^{\nu(|b|-1)}}\int_{R^{\nu(|b|-1)}}e^{i(\varphi_b(x_b,\xi_b)-y_b\cdot\xi_b)}f(y_b)dy_bd\xi_b\tag 5.17
$$
for $f\in \HH_b=L^2(R^{\nu(|b|-1)})$ as an oscillatory integral (see e.g. [KK]).
Wave operators $W_b^\pm$ are now defined by
$$
W_b^\pm=\text{s-}\lim_{t\to\pm\infty}e^{itH}J_be^{-itH_b}P_b.\tag 5.18
$$
We note that this definition depends on $\theta,d,R_0$, but applying stationary phase method to $e^{-itT_b}$ in $e^{-itH_b}=e^{-itT_b}\otimes e^{-itH^b}$ on the RHS we see that the dependence disappears in the limit $t\to\pm\infty$ by the remark made just before Theorem 5.2. Further the asymptotic behavior seen by the stationary phase method tells that the inclusion (5.1) holds:
$$
\RR(W_b^\pm)\subset S_b^0,\tag 5.19
$$
if the limits (5.18) exist. The existence of (5.18) follows from Theorem 5.2-iii), the asymptotic behavior of $e^{-itT_b}$ and Assumptions 1.1-1.2 by noting the relations
$$
\split
&(HJ_b-J_bH_b)e^{-itH_b}P_bf(x)=((T_b+I_b(x_b,x^b))J_b-J_bT_b)e^{-itH_b}P_bf(x),\\
&((T_b+I_b(x_b,0))J_b-J_bT_b)g(x_b)\\
&\hskip20pt =(2\pi)^{-\nu(|b|-1)}\int_{R^{\nu(|b|-1)}}\int_{R^{\nu(|b|-1)}}e^{i(\varphi_b(x_b,\xi_b)-y_b\cdot\xi_b)}a_b(x_b,\xi_b)g(y_b)dy_bd\xi_b
\endsplit\tag 5.20
$$
and the fact that $\text{s-}\lim_{M\to\infty}P_b^M=P_b$. Thus to prove the reverse inclusion (5.2)
$$
S_b^0\subset\RR(W_b^\pm),\tag 5.21
$$
it suffices to prove that
$$
f\in S_b^0\ominus \RR(W_b^\pm)\tag 5.22
$$
implies
$$
f=0.\tag 5.23
$$
To see this we consider the case $t\to+\infty$ and the quantity
$$
(I-e^{isH}J_be^{-isH_b}J_b^{-1})e^{-itH}f=(J_b-e^{isH}J_be^{-isH_b})J_b^{-1}e^{-itH}f\tag 5.24
$$
for $f\in S_b^0(\Delta)$, $\Delta\Subset R^1-\TT$, and $t,s\ge 0$ and use Enss method. Here the existence of $J_b^{-1}$ follows from Theorem 3.3 of [K4] by taking $R_0>0$ in (5.12) large enough (with a slight adaptation to the present case for phases and symbols satisfying the estimates in Theorem 5.2). (5.24) equals
$$
-i\int_0^s e^{iuH}(HJ_b-J_bH_b)e^{-iuH_b}J_b^{-1}du\ e^{-itH}f.\tag 5.25
$$
By Definition 2.1-i)-(2.4) of $S_b^0(\Delta)$, we approximate $f$ by $h \in S_b^{0\sigma}(\Delta)$ for some small $\sigma>0$ with an arbitrarily small error $\delta>0$ so that $\Vert f-h\Vert<\delta$. Then we have for any sufficiently large $R>0$
$$
\limsup_{t\to\infty}\left\Vert e^{-itH}h-\prod_{\alpha\not\le b}F(|x_\alpha|\ge \sigma t)F(|x^b|\le R)e^{-itH}h\right\Vert<\delta.\tag 5.26
$$
Proposition 2.2-iii) and $h \in S_b^{0\sigma}(\Delta)$ yield that for some sequence $t_m\to\infty$ (as $m\to\infty$)
$$
\Vert(\varphi(x_b/t_m)-\varphi(v_b))e^{-it_mH}h\Vert\to0 \q \text{as}\q m\to\infty
$$
for any $\varphi\in C_0^\infty(R^{\nu(|b|-1)})$. Replacing $t$ and $f$ in (5.25) by $t_m$ and $h$, we can therefore insert or remove the factor
$$
\Phi=\prod_{k=1}^{k_b}Q_k\tilde F(|p_{bk}|\ge \sigma')\tilde F(|p_b|\le S)\tilde F(|z_{bk}|/t\ge \sigma')\tilde F(|x^b|\le R) \tag 5.27
$$
to or from the left of $e^{-it_mH}h$ in (5.25) anytime with an error $\delta>0$. Here $p_{bk}=\frac{1}{i}\frac{\partial}{\partial z_{bk}}$, $\sigma'>0$ is a small number with $\sigma'<\sigma$, $\tilde F(|p_b|\le S)$ comes from $E_H(\Delta)$ in $h=E_H(\Delta)h$, and $\tilde F(\tau\le S)$ is a smooth characteristic function of the set $\{\tau\in R^1|\ \tau\le S\}$ with a slope independent of $S$, and $Q_k$ is a pseudodifferential operator
$$
Q_kg(x_b)=(2\pi)^{-\nu(|b|-1)}\int_{R^{\nu(|b|-1)}}\int_{R^{\nu(|b|-1)}}
e^{(x_b\cdot \xi_b-y_b\cdot\xi_b)}
q_k(z_{bk},\zeta_{bk})g(y_b)dy_bd\xi_b\tag 5.28
$$
with symbol $q_k(z_{bk},\zeta_{bk})$ satisfying
$$
\left.\matrix
\format &\l \\
&|\partial_{z_{bk}}^\beta\partial_{\zeta_{bk}}^\gamma q_k(z_{bk},\zeta_{bk})|
\le C_{\beta\gamma}\langle z_{bk}\rangle^{-|\beta|}
\langle \zeta_{bk}\rangle^{-|\gamma|},\\
&q_k(z_{bk},\zeta_{bk})=0\q \text{for}\q \cos(z_{bk},\zeta_{bk})\le \theta\ \text{or} \ |z_{bk}|\le R_0.
\endmatrix \right.
\tag 5.29
$$
The order of products in (5.25) of factors in (5.27) and $J_b^{-1}$ may be arbitrary because these factors are mutually commutative asymptotically as $t\to\infty$ by virtue of (5.26).
We note that $d>0$ in the definition of $J_b=J_{b,\theta,d,R_0}$ can be taken smaller than $\sigma' >0$ beforehand since $W_b^+$ is independent of $d>0$ as mentioned. Thus we can assume the following in addition to (5.29):
$$
q(z_{bk},\zeta_{bk})=0\q \text{for}\q |\zeta_{bk}|\le d.\tag 5.30
$$

We now insert the decomposition (1.21) to the left of $e^{-it_mH}h$ in (5.25) with noting $(I-P^{M_1^m})f=f$ by $f\in \HH_c(H)$. Then by $\Vert(I-P^{M_1^m})h-h\Vert<2\delta$ and by inserting the factor (5.27) to the left of $e^{-it_mH}h$ after the insertion of (1.21), we have
$$
\limsup_{m\to\infty}\Vert(I-P_b^{M_{|b|}^m})\Phi e^{-it_mH}h\Vert<3\delta.\tag 5.31
$$
By the factor $\tilde F(|x^b|\le R)$ in (5.27) and $E_H(\Delta)$ in $h=E_H(\Delta)h$, $P_b^{M_{|b|}^m}$ in (5.31) converges to $P_b$ as $m\to\infty$ in operator norm in the expression (5.31).
It thus suffices to consider the quantity
$$
\int_0^s e^{iuH}((T_b+I_b(x_b,x^b))J_b-J_bT_b)e^{-iuH_b}J_b^{-1}du\ P_b^{M_{|b|}^{m_0}}\Phi e^{-it_mH}h\tag 5.32
$$
for some large but fixed $m_0$ with an error $\delta>0$.
Since $P_b^{M_{|b|}^{m_0}}=\sum_{j=1}^{M_{|b|}^{m_0}} P_{b,E_j}$ $(0\le M_{|b|}^{m_0} <\infty)$ with $P_{b,E_j}$ being one dimensional eigenprojection of $H^b$ corresponding to eigenvalue $E_j$, (5.32) is reduced to considering
$$
\int_0^s e^{-iuE_j} e^{iuH}((T_b+I_b(x_b,x^b))J_b-J_bT_b)P_{b,E_j}e^{-iuT_b}J_b^{-1}\Phi du\  e^{-it_mH}h.\tag 5.33
$$
By Assumptions 1.1-1.2, the factor $P_{b,E_j}$ bounds the variable $x^b$ and yields a short-range error of order $O(\left(\min\langle z_{bk}\rangle\right)^{-1-\ep})$ on the left of $e^{-iuT_b}$ when we replace $I_b(x_b,x^b)$ by $I_b(x_b,0)$, and we have that (5.33) equals
$$
\split
\int_0^s e^{-iuE_j} e^{iuH}P_{b,E_j}O(\langle x^b\rangle)((T_b+I_b(x_b,0))J_b-J_bT_b+&O(\left(\min\langle z_{bk}\rangle\right)^{-1-\ep}))\\
&\times e^{-iuT_b}J_b^{-1}\Phi du\  e^{-it_mH}h,
\endsplit
\tag 5.34
$$
where $O(\langle x^b\rangle)$ is an operator such that $\langle x^b\rangle^{-1}O(\langle x^b\rangle)$ is bounded.
Using (5.20) and the estimate (5.16) in Theorem 5.2-iii) and applying the propagation estimates in Lemma 3.3-ii) of [IK] (again with a slight adaptation to the present case), we now get the estimate:
$$
\Vert ((T_b+I_b(x_b,0))J_b-J_bT_b+O(\left(\min\langle z_{bk}\rangle\right)^{-1-\ep}))
e^{-iuT_b}J_b^{-1}\Phi\left(\min\langle z_{bk}\rangle\right)^{\ep/2}\Vert \le C\langle u\rangle^{-1-\ep/2}
\tag 5.35
$$
for some constant $C>0$ independent of $u\ge 0$. On the other hand (5.26) yields
 that
$$
\Vert\left(\min\langle z_{bk}\rangle\right)^{-\ep/2}e^{-it_mH}h\Vert
$$
is asymptotically less than $2\delta$ as $m\to\infty$.
This and (5.35) prove that the norm of (5.32) is asymptotically less than a constant times $\delta$ as $m\to\infty$.

Returning to (5.24) we have proved that
$$
\split
&\limsup_{m\to\infty}\sup_{s\ge 0}\Vert (I-e^{isH}J_b e^{-isH_b}J_b^{-1})e^{-it_mH} f\Vert\\
&\q \approx_\delta \limsup_{m\to\infty}\sup_{s\ge 0}\Vert (I-e^{isH}J_b e^{-isH_b}J_b^{-1})P_b^{M_{|b|}^m}e^{-it_mH} f\Vert\le C \delta,
\endsplit
\tag 5.36
$$
where $a\approx_\delta b$ means that $|a-b|\le C\delta$ for some constant $C>0$.
Since wave operator $W_b^+=\text{s-}\lim_{s\to\infty}e^{isH}J_be^{-isH_b}P_b$ exists, (5.36) yields
$$
\limsup_{m\to\infty}\Vert (I-W_b^+J_b^{-1})P_b^{M_{|b|}^m}e^{-it_mH}f\Vert \le C\delta.\tag 5.37
$$
By the arguments above deriving (5.31) we can remove $P_b^{M_{|b|}^m}$ and get
$$
\limsup_{m\to\infty}\Vert (I-W_b^+ J_b^{-1})e^{-it_mH}f\Vert \le C\delta.\tag 5.38
$$
Since we assumed (5.22), $f$ is orthogonal to $\RR(W_b^+)$. Thus taking the inner product of the vector inside the norm in (5.38) with $e^{-it_mH}f$, we have
$$
\Vert f\Vert^2=\lim_{m\to\infty}|(e^{-it_mH}f, e^{-it_mH}f)|
=\lim_{m\to\infty}|(e^{-it_mH}f,(I-W_b^+ J_b^{-1})e^{-it_mH}f)|
\le C\delta\Vert f\Vert.
$$
As $\delta>0$ is arbitrary, this gives $f=0$, proving (5.23). The proof of (4.116) is complete.

\BP

\enddocument